\documentclass[reqno,english]{amsart}
\usepackage{geometry}

\usepackage[maxbibnames=99, giveninits=true, style=numeric]{biblatex} 
\renewbibmacro{in:}{}
\usepackage{tikz}
\usetikzlibrary{angles,quotes}

\usepackage{comment}

\usepackage{amssymb}
\usepackage{amsthm}
\usepackage{amsmath}
\usepackage{bbm}
\usepackage{bm}
\usepackage[bottom]{footmisc}
\usepackage{braket}
\usepackage{enumitem}
\usepackage{esint}
\usepackage{extarrows}
\usepackage{graphicx}    
\usepackage{hyperref}
\usepackage{ifthen}
\usepackage{mathrsfs}
\usepackage{mathtools}
\usepackage[version=3]{mhchem}
\usepackage{microtype}
\usepackage{stmaryrd}
\usepackage{tikz-cd}
\usepackage{xurl} 

\usepackage[T1]{fontenc}
\usepackage{mlmodern}

\newcommand{\abs}[1]{\left\lvert#1\right\rvert}

\newcommand{\brac}[1]{\left[ #1 \right]}
\newcommand{\braces}[1]{\left\{ #1 \right\}}

\newcommand{\inner}[1]{\langle #1\rangle}
\newcommand{\paren}[1]{\left( #1 \right)}

\newcommand{\sbig}[1]{\vcenter{\hbox{\scalebox{1.10}{\ensuremath{#1}}}}}
\newcommand{\sparen}[1]{\sbig( #1 \sbig)}

\newcommand{\eps}{\varepsilon}

\newcommand{\loc}{\text{loc}}
\newcommand{\p}{\partial}

\newcommand{\sm}{\setminus}

\newcommand{\xto}[1]{\xrightarrow{#1}}

\newcommand{\pfrac}[2]{\left(\frac{#1}{#2}\right)}

\newcommand{\wt}{\widetilde}

\newcommand{\R}{\mathbb{R}}

\renewcommand{\phi}{\varphi}

\let\originalleft\left
\let\originalright\right
\renewcommand{\left}{\mathopen{}\mathclose\bgroup\originalleft}
\renewcommand{\right}{\aftergroup\egroup\originalright}

    \let\div\odiv

\DeclareMathOperator{\div}{div}

\DeclareMathOperator{\Tr}{Tr}

\makeatletter
\def\blfootnote{\xdef\@thefnmark{}\@footnotetext}
\makeatother

\newlist{aenumerate}{enumerate}{50}
\setlist[aenumerate]{label=(\alph*)}

\newlist{penumerate}{enumerate}{50}
\setlist[penumerate]{label=(\arabic*)}

\newlist{renumerate}{enumerate}{50}
\setlist[renumerate]{label=(\roman*)}

\theoremstyle{plain}
\newtheorem{thm}{Theorem}[section]
\newtheorem{cor}[thm]{Corollary}

\newtheorem{lemma}[thm]{Lemma}

\theoremstyle{definition}

\AtEndEnvironment{example}{\null\hfill\qedsymbol}

\theoremstyle{remark}
\newtheorem{remark}{Remark}
\AtEndEnvironment{remark}

\allowdisplaybreaks[1]

\usepackage{mathtools}
\mathtoolsset{showonlyrefs}

\title[Derivation of the Quantum Landau Operator]{From the quantum Boltzmann operator to the Quantum Landau operator}

\author{Maria Pia Gualdani}
\address{The University of Texas at Austin}
\email{gualdani@math.utexas.edu}

\author{Nata\v sa Pavlovi\' c}
\address{The University of Texas at Austin}
\email{natasa@math.utexas.edu}

\author{Justin Toyota}
\address{The University of Texas at Austin}
\email{jtoyota@utexas.edu}

\author{Dominic Wynter}
\address{The University of Texas at Austin}
\email{dominic.wynter@austin.utexas.edu}

\begin{document}

\begin{abstract}
    In this manuscript we derive the quantum Landau operator as the weak-coupling limit of the quantum Boltzmann operator (also known as the Uehling-Uhlenbeck operator). We consider both Fermi-Dirac and Bose-Einstein statistics. Our approach is inspired by the work by Benedetto and Pulvirenti \cite{BenedettoPulvirenti07}, where the classical Landau operator was derived from the quantum Boltzmann operator. To capture the ternary term in the quantum Landau operator, we introduce a new two-parameter scaling that preserves the quantum effects in the limit. Furthermore, we provide an explicit rate of convergence that depends on the regularity of the interaction potential.   
\end{abstract}

\maketitle

\setcounter{tocdepth}{2}
 
\let\oldtocsection=\tocsection
 
\let\oldtocsubsection=\tocsubsection
  
\renewcommand{\tocsection}[2]{\hspace{0em}\oldtocsection{#1}{#2}}
\renewcommand{\tocsubsection}[2]{\hspace{1em}\oldtocsubsection{#1}{#2}}
\renewcommand{\tocsubsubsection}[2]{\hspace{2em}\oldtocsubsubsection{#1}{#2}}

\setcounter{tocdepth}{2}
\tableofcontents

\section{Introduction}

The focus of this paper is a derivation of the quantum Landau operator from the quantum Boltzmann operator. Let us start by introducing these operators and the equations containing them. 

The quantum Boltzmann equation (also known as the Uehling-Uhlenbeck or Boltzmann-Nordheim equation) 
was introduced by Nordheim \cite{Nordheim28} and Uehling and Uhlenbeck \cite{UU33} to model a dilute gas of weakly-interacting bosons or fermions, generalizing the classical Boltzmann equation derived by Maxwell \cite{Maxwell67} and Boltzmann \cite{Boltzmann72}. It reads as 
\begin{equation} \label{eq:UU_formal_def_gen}
\begin{split}
    \p_{t}f + v \cdot \nabla_x f 
    &= \; \frac{1}{8\pi^{2}\hbar^{2}}\int_{\R^{3}} \int_{\R^{3}}
    \brac{\hat{V}(k) +\theta \hat{V}\paren{k + \frac{u}{\hbar}}}^{2}
    \delta(\hbar |k|^{2} + k\cdot u) 
    \\
    &\qquad\qquad\qquad\qquad\quad\cdot
    \Big\{(1 +\theta (2\pi\hbar)^3 f)(1 +\theta (2\pi\hbar)^3 f_{1})f'\!f_{1}'
    \\
    &\qquad\qquad\qquad\qquad\qquad\quad
    - (1 +\theta (2\pi\hbar)^3 f')(1 +\theta (2\pi\hbar)^3 f'_{1})f\! f_{1}\Big\}\,dk\,dv_{1}
\end{split}
\end{equation} 
for an interaction potential $V$ and a quantum parameter $\hbar>0$, where we have written $u = v-v_{1}$ and taken $\theta=\pm1$, with $\theta=1$ modeling Bose-Einstein particles and $\theta=-1$ modeling Fermi-Dirac particles. We have also used the standard notation
\begin{align*}
    f &= f(v) & f' &= f(v') &
    f_{1} &= f(v_{1}) & f'_{1} &= f(v'_{1}),
\end{align*}
where
\begin{align}
    \label{eq:postcollisional_velocities}
    v' &= v + \hbar k,
    &
    v'_{1} &= v_{1} - \hbar k
\end{align}
are the post-collisional velocities. We compare to the classical Boltzmann equation, which can be written as
\begin{align}
    \label{eq:classical_Boltzmann}
    \p_t f + v\cdot\nabla_x f 
    =
    \int_{\R^3}\!\int_{S^2}\! B\left(v-v_1,\omega\right)[f' f_1' - f\! f_1]\, d\omega\, d v_1
\end{align}
for a given collision rate $B(v-v_1,\omega)$ and post-collisional velocities
\begin{align*}
    v' = v - \langle \omega, v-v_1\rangle \omega,\qquad
    v_1' = v_1 + \langle \omega, v-v_1\rangle \omega.
\end{align*}
Whereas the classical Boltzmann operator in \eqref{eq:classical_Boltzmann} is purely binary, the quantum Boltzmann operator in \eqref{eq:UU_formal_def_gen} has both binary and ternary parts, representing the effect of Fermi blocking or Bose enhancement. 

We now summarize the current analytical results for the quantum Boltzmann equation \eqref{eq:UU_formal_def_gen}, which can be divided into the fermionic ($\theta=-1$) and bosonic ($\theta=1$) cases. The two cases exhibit vastly different properties, reflecting the different physics of fermions compared to bosons. Broadly speaking, the former case reflects the Pauli exclusion principle, while the latter exhibits the condensation properties characteristic of bosons. 

In the fermionic case, the first works exhibiting existence of global weak solutions to \eqref{eq:UU_formal_def_gen} are due to Dolbeault \cite{Dolbeault1994} and Lions \cite{Lions1994}. 
In the spatially homogeneous setting, existence of global smooth solutions and convergence to equilibrium were proved by Escobedo, Mischler and Valle \cite{EscobedoMischlerValle09}. 

In the bosonic case, we first mention the results obtained in the homogeneous case. In that context, most research has focused on isotropic (radial) initial data. Namely, Lu \cite{lu2000} proved global existence of an isotropic mild $L^1$ solution under strong conditions on the collision kernel. Under weaker conditions on the kernel, which include the hard sphere case, Lu subsequently proved the existence of isotropic measure-valued solutions \cite{lu2004isotropic} and then proved weak convergence to their associated (potentially singular) Bose-Einstein equilibria in \cite{lu2005boltzmann}. For non-radial data, local in time strong solutions were constructed by Briant and Einav \cite{BriantEinav16}, and for high-temperature initial data, global strong solutions were constructed by Li and Lu \cite{LiLu19}. On the other hand, Escobedo and Vel\'azquez proved that there exist solutions to \eqref{eq:UU_formal_def_gen} that exhibit $L^\infty$ blowup in finite time \cite{EscobedoVelazquez15} \cite{EscobedoVelazquez14} for initial data with temperature $T<T_{crit}$ below the critical value. This blow-up singularity corresponds to the physical phenomenon of condensation. 

In the inhomogeneous case, global well-posedness of \eqref{eq:UU_formal_def_gen} for Bose-Einstein statistics is known for initial data near equilibrium with temperature above the critical value \cite{OuyangWu22} \cite{zhou2023global}. 

Derivation of the quantum Boltzmann equation 
\eqref{eq:UU_formal_def_gen}
is a challenging mathematical problem that still attracts attention. Early derivation results relied on truncations of series expansions arising from the many-particle dynamics of lattice Fermi gases, such as the Dyson-series methods of Hugenholtz \cite{Hugenholtz83} and Ho and Landau \cite{HoLandau97}. Later, important progress was obtained by Erd\H{o}s, Salmhofer and Yau \cite{ESY04} and Benedetto, Castella, Esposito and Pulvirenti \cite{BCEP04} \cite{BCEP05} \cite{BCEP06} \cite{BCP07} \cite{BCEP08}, using a weak-coupling rescaling of the Wigner transform of a solution to the many-body Schr\"{o}dinger equation. More precisely, in \cite{BCEP04} \cite{BCEP08} the authors prove convergence to the binary Boltzmann hierarchy for factorized initial data satisfying Maxwell-Boltzmann statistics, and in \cite{BCEP05} prove convergence of a truncated Duhamel expansion to the quantum Boltzmann equation, for initial data corresponding to quasi-free states. On the other hand, the work of Erd\H{o}s, Salmhofer and Yau \cite{ESY04}
proceeds by taking limits of the one-particle marginal distributions under the assumption that the solution of the many-body system propagate quasi-freeness in time and proving convergence to the quantum Boltzmann equation (cf. \cite{LukkarinenSpohn09} for a review of these results). These analyses were revisited by X.\ Chen and Guo \cite{ChenGuo15} and X.\ Chen and Holmer \cite{ChenHolmer23} who proved, conditionally on sufficient regularity of solutions to the BBGKY hierarchy, that such solutions will converge to the classical Boltzmann hierarchy rather than the quantum Boltzmann hierarchy. Recently, the quantum Boltzmann operator in \eqref{eq:UU_formal_def_gen} was derived for small timescales for the first time without assuming quasi-freeness by T.\ Chen and Hott \cite{ChenHott23, ChenHott25} in the case of bosons, and by T.\ Chen and Card\'{e}nas \cite{CardenasChen23} in the case of fermions.

A closely related equation is the quantum Landau equation, which we write as
\begin{align}  \label{quantum_Landau eq}
    \p_t f + v\cdot \nabla_x f &= Q_{qL}^{\alpha_{0}}(f),
\end{align}
with the operator $Q_{qL}^{\alpha_{0}}$ defined via
\begin{align}    \label{quantum_Landau}
    Q_{qL}^{\alpha_{0}}(f) := B \operatorname{div}_{v}\int_{\R^{3}}\left[ K(v-v_1)\big(f_1(1+\theta\alpha_{0} f_1)\nabla_v f-f(1+\theta\alpha_{0} f)\nabla_{v_1}f_1\big)\right]\,dv_{1}
\end{align}
for a rate parameter $B>0$, a quantum parameter $\alpha_0 \geq 0$ and $\theta = \pm 1$, and where $K(v-v_1)$ is the matrix defined by
\begin{align}
    \label{angular_projection}
     K(z) := |z|^{-1}(I - \hat{z}\otimes \hat{z}),
     \qquad \hat{z} := \frac{z}{|z|},\,z\neq0.
\end{align}


\noindent The fermionic case ($\theta =-1$) of \eqref{quantum_Landau}, known as the Landau-Fermi-Dirac equation, has been extensively studied, and well-posedness for spatially homogeneous data is known for hard potentials due to Bagland \cite{bagland2004well}, for moderately soft potentials due to Alonso, Bagland, Desvillettes and Lods \cite{alonso2022landau}, and for Coulomb potential due to Golding, Gualdani and Zamponi \cite{GoldingGualdaniZamponi22}. More recently, global well-posedness of strong solutions for the inhomogeneous equation \eqref{quantum_Landau} was established by Golding and Henderson \cite{GoldingHenderson25}. Meanwhile, the bosonic case ($ \theta =1$) of \eqref{quantum_Landau}, known as the Landau-Bose-Einstein equation, finds applications in models of dark matter \cite{Chavanis2021}. However, this case has received relatively little mathematical study. Recently, the first, second and fourth authors \cite{GualdaniPavlovicWynter25} showed existence of smooth solutions to \eqref{quantum_Landau} with $ \theta =1$ with initial data near equilibrium. 

The first connection between the quantum Boltzmann equation \eqref{eq:UU_formal_def_gen} and 
the quantum Landau equation \eqref{quantum_Landau eq} was established by Danielewicz \cite{Danielewicz80}; there the quantum Landau operator \eqref{quantum_Landau} was derived by a formal grazing limit of the quantum Boltzmann operator (the right side of \ref{eq:UU_formal_def_gen}). This brings us to the main theme of this paper, namely a rigorous derivation of the quantum Landau operator from the quantum Boltzmann operator in the bosonic and fermionic regimes by a weak-coupling high-density limit. 
The approach that we take is inspired by 
the works on semiclassical limits of the quantum Boltzmann operator appearing in \eqref{eq:UU_formal_def_gen}
by Benedetto and Pulvirenti \cite{BenedettoPulvirenti07} and He, Lu and Pulvirenti \cite{he2021semi}, where convergence to the operator appearing in the Landau-Fokker-Planck equation was proved by a weak-coupling-type limit. 

However, in order to derive ternary terms in the quantum Landau operator \eqref{quantum_Landau} we introduce a new two parameter family of scalings. Before we present details, we first introduce our notation in the subsection below.

\subsection{Notation}We will use the following notation for Sobolev norms: for any function $f\in L^1_{\mathrm{loc}}(\mathbb{R}^3)$ which is $k$-times weakly differentiable for $k\in\mathbb{Z}_{\ge0}$, we define the norm
\begin{align*}
    \| f\|_{W^{k,p}(\mathbb{R}^3)}
    =
    \begin{cases}
        \left( \sum\limits_{|\alpha|\le k}\|\partial_v^\alpha f\|_{L^p(\mathbb{R}^3)}^p\right)^{1/p},
        &
        p\in[1,\infty)
        \\
        \max\limits_{|\alpha|\le k}\| \partial_v^\alpha f\|_{L^\infty(\mathbb{R}^3)},
        & p = +\infty
    \end{cases}
\end{align*}
and we define $W^{k,p}(\mathbb{R}^3)$ to be the space of functions $f$ such that $\| f\|_{W^{k,p}(\mathbb{R}^3)}<\infty$. 

The interaction potential $V$ appears in our operators via its Fourier transform
\begin{equation*}
    \hat{V}(k) := \int_{\R^{3}}e^{-ik\cdot v}V(v)\,dv.
\end{equation*}
Throughout this paper we will assume that $V$ is radial, which implies $\hat{V}$ is also radial. By abuse of notation, we can thus view $\hat{V}$ as a function on $[0,\infty)$; i.e., we write $\hat{V}(r)$ to denote the value of $\hat{V}$ at any $k\in \R^{3}$ with $|k|=r$. The primary functions of $\hat{V}$ that we will deal with are the moments
\begin{align}
    \label{defn_of_M_a}
    M_a = \int_0^\infty\! (\mu+\mu^a)\hat V(\mu)^2\, d\mu;
\end{align}
note that $V$ being even implies $\hat{V}$ is real-valued, so this is always non-negative. Finally, in discussing the hypotheses of our main theorem, we will use the Japanese bracket notation $\inner{v} := \sqrt{1 + |v|^{2}}$. We will at times use the notation
\begin{equation*}
    \|f\|_{L^{\infty}_{s}} := \|\inner{v}^{s}f\|_{L^{\infty}},
    \qquad s\in \R.
\end{equation*}

\subsection{Main Result}

We now introduce a two parameter family of scalings $(\eps, \alpha)$ in the quantum Boltzmann operator, as follows:
for a given $\eps>0$, $\alpha = \alpha(\eps)\geq 0$ and sufficiently regular radial interaction potential $V$, we define the $Q_{qB}^{\alpha,\varepsilon}$ operator via
\begin{equation} \label{eq:UU_formal_def}
\begin{split}
    Q_{qB}^{\alpha,\varepsilon}(f)(v)
    :=& \frac{1}{8\pi^{2}\eps^{2}}\int_{\R^{3}} \int_{\R^{3}}
    \brac{\hat{V}(k) +\theta \hat{V}\paren{k + \frac{u}{\eps}}}^{2}
    \delta(\eps |k|^{2} + k\cdot u) 
    \\
    &\quad\cdot
    \{(1 +\theta \alpha f)(1 +\theta \alpha f_{1})f'\!f_{1}'
    - (1 +\theta \alpha f')(1 +\theta \alpha f'_{1})f\! f_{1}\}\,dk\,dv_{1}
\end{split}
\end{equation}
for sufficiently regular functions $f$ (for precise conditions on $V$ and $f$, see Theorem \ref{mainthm}),
where we have written $u = v-v_{1}$, and where $\theta=\pm1$.
We interpret $Q_{qB}^{\alpha,\varepsilon}(f)$ as acting on test functions $\psi(v)\in W^{3,\infty}(\R^{3})$ as follows 
(details can be found in Subsection \ref{sect:expansion_of_uu}):
\begin{equation}
\begin{aligned} \label{eq:UU_weak_form}
    &\int_{\R^3}  Q_{qB}^{\alpha,\varepsilon}(f)(v)\psi(v)\,dv
    \\
    &\quad :=
    \frac{1}{4\pi^2\eps^4}\int_{\R^3\times\R^3} f\! f_1\int_{S^2_-} |\hat k\cdot u|
    \hat{V}\left(\frac{|\hat k\cdot u|}{\eps}\right)\left[ \hat{V}\left(\frac{|\hat k\cdot u|}{\eps}\right)
    +\theta \hat{V}\left(\frac{\sqrt{u^2-|\hat k\cdot u|}}{\eps}\right)\right] 
    \\
    &\qquad\qquad\qquad\qquad\qquad\qquad\qquad\qquad\cdot[\psi(v')-\psi(v)](1+\theta\alpha[f'+f_1'])\,d\hat{k}\,dv_{1}\,dv,
\end{aligned}
\end{equation}
where the half-sphere $S^2_{-}$ is defined by
\begin{equation}
\label{eq:half_sphere_definition}
    S_{-}^{2} := S_{-}^{2}(u)
    = \{\hat{k}\in S^{2} : \hat{k}\cdot u \leq 0\}.
\end{equation}
In order to state our main results we also need a weak formulation of the quantum Landau operator \eqref{quantum_Landau}, which we establish in Lemma \ref{weak_formulation_quantum_Landau}, and record now: 
\begin{align} 
        &\int_{\R^3}\! 
        Q_{qL}^{\alpha_0}(f)(v)\psi(v)\,dv 
        \label{intro-qL-weak} \\
        &= \int_{\R^3\times\R^3}\!
        \left[ - \frac{2B (v-v_1)}{|v-v_1|^3}\cdot(\nabla_v \psi -\nabla_{v_1}\psi_1) + B\Tr(K(v-v_1)D^2 \psi(v))\right] f\! f_1(1+\theta\alpha_0 f_1)\,dv_{1}\,dv, \nonumber
\end{align}
for a rate parameter $B>0$ and quantum parameter $\alpha_0\ge0$, and with the matrix $K$ as defined in \eqref{angular_projection}. 
Our main results are summarized in the following theorem.

\begin{thm}  \label{mainthm}
Let $\eta \in (0,1)$ and $0\not\equiv V\in L^2(\mathbb{R}^3)$ be a radial potential which satisfies 
\begin{equation} \label{eq:Vcond}
    M_{5+2\eta} := \int_0^\infty(\mu+\mu^{5+2\eta})\hat{V}(\mu)^2\,d\mu<\infty.
\end{equation} 
Let $\alpha=\alpha(\eps)\ge0$ be such that 
\begin{equation*}
    \lim_{\eps\to0}\alpha(\eps)=\alpha_0\in[0,\infty).
\end{equation*}
Then for any $\psi\in W^{3,\infty}(\R^3)$ and $f$ such that $\langle v\rangle^N f\in W^{2,\infty}(\mathbb{R}^3)$ for some $N>3$, we have the estimate
\begin{multline} \label{eq:convergence_rate}
    \left|
        \int_{\R^3}\!
        \left( Q_{qB}^{\alpha,\eps}(f)-Q_{qL}^{\alpha_0}(f)
        \right)\psi(v)\, dv
    \right| \\
    \lesssim_{N,\eta} (\eps^\eta + |\alpha-\alpha_0|)
    M_{5+2\eta}
    (1 + \alpha_0 + \alpha)
    \|\psi\|_{W^{3,\infty}}
    \paren{\|\inner{v}^{N}f\|_{W^{1,\infty}}^2 + \|\inner{v}^{N}f\|_{W^{2,\infty}}^3},
\end{multline}
where $\int_{\R^{3}} Q_{qB}^{\alpha,\varepsilon}(f)\psi(v)\,dv$ and $\int_{\R^{3}} Q_{qL}^{\alpha_0}(f)\psi(v)\,dv$ are defined in (\ref{eq:UU_weak_form}) and (\ref{intro-qL-weak}), respectively, with constant
\begin{equation} \label{eq:Bdef}
    B = \frac{1}{8\pi}\int_{0}^{\infty}\mu^{3}\hat{V}(\mu)^{2}\,d\mu
\end{equation}
in $Q_{qL}^{\alpha_0}$. Consequently,
    \begin{align} \label{eq:main limit}
        \int_{\R^{3}} Q_{qB}^{\alpha,\varepsilon}(f)\psi(v)\,dv \to \int_{\R^{3}} Q_{qL}^{\alpha_0}(f)\psi(v)\,dv ,\qquad \textrm{as}\; \eps\to 0. 
    \end{align}
\end{thm}

\begin{remark}
For fermions, when $\alpha_0\ne0$, it seems to us that the limit obtained in the above theorem could be understood as a high-density limit, due to the implicit rescaling of velocity $(v,v_*)\mapsto (\eps v,\eps v_*)$ in \eqref{eq:UU_weak_form}, which can be viewed as taking data $f$ with mass $M\sim\eps^{-3}$ and energy $E\sim\eps^{-5}$ as $\eps \rightarrow 0$. In addition, the limit has some aspects of a grazing collision limit, which is used to connect the classical Boltzmann operator with the classical Landau operator. More precisely, the collision term weak formulation \eqref{eq:UU_weak_form} formally implies that collisions concentrate on the set where $|u|\cos\vartheta=|\hat k\cdot u| \lesssim\eps$, so that collisions with small relative velocity $u$, or with angle of incidence $\vartheta\approx\pi/2$, dominate. This can be compared to foundational grazing limit results such as \cite{Desvillettes92}, which in our parametrization can be interpreted as taking collisions concentrated on the set where $\abs{\vartheta-\pi/2}\lesssim\eps$. Our result, in contrast, involves a rescaled interaction potential, taking both the deflection angle $\vartheta'=2\vartheta-\pi$ and deflection velocity $|u|$ to zero as $|u|\sin(2\vartheta)\to0$, with the result that we derive the Landau equation with Coulomb potentials, regardless of the original interaction potential.
\end{remark}

\begin{remark}
We expect that a similar result might be obtained when $\langle v\rangle^N f\in L^\infty$ for $N>3$ with a potentially weaker rate of convergence, though we will not pursue this path in the current work. 
\end{remark}

As a special case of the above theorem, when $\alpha = (2\pi\eps)^3$ we can recover the result of \cite{BenedettoPulvirenti07}, where the classical Landau equation was derived as a limit. More precisely, we obtain: 

\begin{cor}
    Let the potential $ V$ and the function $f$ satisfy the same assumptions as in Theorem \ref{mainthm}, and let $
    \psi\in W^{3,\infty}(\mathbb{R}^3)$. Then we have the limit
    \begin{align*}
        \int_{\R^{3}} Q_{qB}^{(2\pi\eps)^3,\eps}(f)\psi(v)\,dv
        \to\int_{\R^{3}}Q_{qL}^{0}(f,f) \psi(v)\,dv,
        \qquad\eps\to0,
    \end{align*}  
    with $B$ given by \eqref{eq:Bdef}, and 
    $Q_{qL}^{0}$ given by 
    \eqref{quantum_Landau} with $\alpha_0 = 0$, corresponding to the classical Landau operator. 
 \end{cor}

Our proof of Theorem \eqref{mainthm} is inspired by the work of Benedetto-Pulvirenti \cite{BenedettoPulvirenti07}. However, in order to capture ternary contributions, in addition to expanding the test function $\psi$ into Taylor series as was done in \cite{BenedettoPulvirenti07}, we further expand the probability distribution $f$ itself. The proof is organized as follows: we start by deriving 
the weak form \eqref{intro-qL-weak} of the quantum Landau operator $Q_{qL}^{\alpha}$ in Subsection \ref{subsec-qL-weak} and the weak form \eqref{eq:UU_weak_form} of the quantum Boltzmann operator $Q_{qB}^{\alpha,\varepsilon}$ in Subsection \ref{sect:expansion_of_uu}. From there, we divide $Q_{qB}^{\alpha,\varepsilon}$ into four terms: 
\begin{itemize}
\item the \textit{binary main term} $Q_{B}^{\eps}$ (which is a collision operator for Maxwell-Boltzmann statistics, see Benedetto-Pulvirenti \cite[p.\ 910]{BenedettoPulvirenti07}), 
\item the \textit{ternary main term} $Q_{T}^{\alpha,\eps}$, and 
\item two \textit{cross terms} $R_{B}^{\eps}$ and $R_{T}^{\alpha,\eps}$.
\end{itemize}
The binary and ternary main terms will be further decomposed into terms of differing asymptotics in Subsections \ref{sect:binary_main} and \ref{sect:ternary_main}, respectively. We then show that a combination of these terms converge to $Q_{qL}^{\alpha_{0}}$, while all of the remaining terms converge to 0 (Subsection \ref{sect:limits_of_main_terms}). Finally, in Subsection \ref{sect:cross_terms}, we show that the cross terms $R_{B}^{\eps}$ and $R_{T}^{\alpha,\eps}$ converge to 0. These limit results are combined in Subsection \ref{subsect:conclusion} to complete the proof of Theorem \ref{mainthm}.

Upon the completion of this work, we learned that a grazing limit for the collision operators in (four and three) wave kinetic equations\footnote{In particular, the collision operator in the four wave kinetic equation corresponds to the ternary part of the quantum Boltzmann operator.} (which appear in the context of wave turbulence) were established by Duong and He \cite{DuongHe2025}, though under a different scaling than the weak-coupling limit considered here.

\subsection{Structure of the Paper}

Section \ref{sect:lemmas} contains a number of calculations that will be used throughout the paper, while the proof of the main theorem (Theorem \ref{mainthm}) is the content of Section \ref{sect:proof_of_main}.

\subsection*{Acknowledgements} MPG is partially supported by the NSF grants DMS-2206677 and 
DMS-2511625. NP is partially supported by the NSF grants DMS-1840314, DMS-2052789 and DMS-2511517. 
JT is partially supported by the NSF grants DMS-1840314, DMS-2052789 and DMS-2511517 through NP. 
This material is based upon work supported by the National Science Foundation under Grant No. DMW-2424139, while the authors were in residence at the Simons Laufer Mathematical Sciences Institute in Berkeley, California during the Fall 2025 semester. 
The authors would like to thank the Simons Laufer Mathematical Institute for their kind hospitality.
Also, we would like to thank Esteban C\'{a}rdenas and Michael Hott for inspiring conversations.

\section{Preliminaries} \label{sect:lemmas}

In this section we introduce preparatory lemmas that we will use in Section \ref{sect:proof_of_main}.

\subsection{Functional Estimate} \label{subsect:functional_estimate}

We start with lemma which shows that the operator $Q_{qB}^{\alpha,\eps}(f)$ is well-defined for slowly-decaying potentials $V$, and particularly that the integral \eqref{eq:UU_formal_def} remains well-defined when $\hat V$ and $f$ only defined almost everywhere.

\begin{lemma}
    \label{lem:QUU_well-defined}
    For any radial potential $V:\R^3\to\R$ such that $$\int_0^\infty\mu\hat V(\mu)^2\, d\mu<\infty,$$
    and any $f$ such that $\inner{v}^{N}f\in L^{\infty}(\R^3)$ for some $N>3$, we have the bound
    \begin{align*}
        \|Q_{qB}^{\alpha,\eps}
        (f)\|_{L^1(\R^3)}
        \lesssim_{N} \eps^{-2}\|f\|_{L^\infty_N}^2\left(1+\alpha\| f\|_{L^\infty}\right)\int_0^\infty\mu\hat V(\mu)^2\,d\mu.
    \end{align*}
\end{lemma}
We note that a similar estimate was proved in \cite{HeLuPulvirentiZhou24}
for a more restrictive class of functions, such that $f$ satisfies the Pauli exclusion principle when Fermi-Dirac statistics are considered. Since the arguments are standard, the proof of the lemma is deferred to the end of Appendix \ref{subsect:delta_identity}.

\subsection{Asymptotic Lemmas} \label{subsect:asymptotic_lemmas}

The following lemma gives formulas for integrals involving expressions of the form  
\begin{equation} \label{eq:def_of_b_a}
    b_{a}(\mu) := \mu^{a}\hat{V}(\mu)^{2},
    \qquad a,\mu\in [0,\infty),
\end{equation}
which will be seen repeatedly in the computations of Subsection \ref{sect:limits_of_main_terms}. These formulas will involve the functions
\begin{align} \label{eq:B_eps}
    B_{a,\eps}(u) := 
    \frac{1}{8\pi}\int_{0}^{|u|/\eps}b_{a}(\mu)\,d\mu = 
    \frac{1}{8\pi}\int_{0}^{|u|/\eps}\mu^{a}\hat{V}(\mu)^{2}\,d\mu,
    \quad\quad u\in \R^{3}
\end{align}
as well as the moments $M_{a}$ defined in \eqref{defn_of_M_a}.

\begin{lemma} \label{lemma:b_a}
Let $u\in \R^{3}\sm \{0\}$ and $\eps>0$. Defining $S_{-}^{2} := \{\hat{k}\in S^{2} : \hat{k}\cdot u \leq 0\}$ and using the notation above, we have:
\begin{aenumerate}
    \item For all $a\geq 0$,
    \begin{equation} \label{eq:b_2}
        \int_{S_{-}^{2}}b_{a}\pfrac{\hat{k}\cdot u}{\eps}\hat{k}\,d\hat{k}
        = -\frac{16\pi^{2}\eps^{2}B_{a+1,\eps}(u)}{|u|^{3}}\,u.
    \end{equation}

    \item For all $a\geq 0$ and $\eta < 2$,
    \begin{equation} \label{eq:b_3}
        \int_{S_{-}^{2}}b_{a}\pfrac{\hat{k}\cdot u}{\eps}\hat{k}\otimes \hat{k}\,d\hat{k}
        = 
        8\pi^{2}\eps(B_{a,\eps}(u)K(u) + R_{\eps}(u)),
    \end{equation}
    where $K(\cdot)$ is defined in \eqref{angular_projection},
    and
    \begin{equation} \label{eq:bound_on_R}
        \|R_{\eps}(u)\| \lesssim M_{a+\eta}\frac{\eps^{\eta}}{|u|^{1+\eta}};
    \end{equation}
    here $\|\cdot\|$ represents any norm on the space of $3\times 3$ real matrices.

    \item For all $a\geq 1$ and $\eta\leq 1$,
    \begin{equation} \label{eq:b_a}
        \int_{S_{-}^{2}}b_{a}\pfrac{\hat{k}\cdot u}{\eps}\,d\hat{k}
        \lesssim M_{a-1+\eta}\pfrac{\eps}{|u|}^{\eta}.
    \end{equation}
\end{aenumerate}
\end{lemma}

The statement of the above lemma is a slight generalization of the statement of Lemma 3.1 in \cite{BenedettoPulvirenti07}. Namely, part ($a$) above allows all $a \geq 0$ rather than only $a=2$, which is the case considered in \cite{BenedettoPulvirenti07}. Likewise, the statement ($b$) above is also valid for all $a \geq 0$, while \cite{BenedettoPulvirenti07} covers only the case when $a=3$.
For these reasons and in order to make the paper self-contained we present the proof of Lemma \ref{lemma:b_a} here. 

\begin{proof}
We begin by proving ($a$). We parametrize $\hat k\in S_{-}^{2}$ via the decomposition
\begin{equation} \label{eq:k_decomp}
    \hat{k} =   \sqrt{1 - \lambda^{2}}\,\xi -\lambda\hat{u},
\end{equation}
where $\hat{u} = u/|u|$ and 
\begin{gather*}
    \xi \in S^{1}(u) := \{v\in S^{2} : v \perp u\}, \\
    \lambda = -\hat{k}\cdot \hat{u} \in (0,1].
\end{gather*}
Then, we can compute that
\begin{align*}
    \int_{S_{-}^{2}}b_{a}\pfrac{\hat{k}\cdot u}{\eps}\hat{k}\,d\hat{k}
    &= \int_{0}^{1}\int_{S^{1}(u)}
    b_{a}\left(-\frac{\lambda|u|}{\eps}\right)\paren{\sqrt{1-\lambda^{2}}\,\xi - \lambda\hat{u}}\,d\xi\,d\lambda \\
    &= \int_{0}^{1}b_{a}\pfrac{\lambda|u|}{\eps}
    \paren{\sqrt{1-\lambda^{2}}\int_{S^{1}(u)}\xi\,d\xi 
    - \lambda\hat{u}\int_{S^{1}(u)}\,d\xi}\,d\lambda \\
    &= -2\pi\int_{0}^{1}\lambda b_{a}\pfrac{\lambda|u|}{\eps}\hat{u}\,d\lambda.
\end{align*}
Using the change of variables $\lambda \mapsto \mu = \lambda|u|/\eps$, this becomes
\begin{align*}
    -2\pi\int_{0}^{1}\lambda b_{a}\pfrac{\lambda|u|}{\eps}\hat{u}\,d\lambda
    &= -\frac{2\pi\eps^{2}}{|u|^{2}}\paren{\int_{0}^{|u|/\eps}\mu b_{a}(\mu)\,d\mu}\hat{u} \\
    &= -\frac{2\pi\eps^{2}}{|u|^{2}}\paren{\int_{0}^{|u|/\eps}\mu^{a+1}\hat{V}(\mu)^{2}\,d\mu}\hat{u} \\
    &= -\frac{16\pi^{2}\eps^{2}B_{a+1,\eps}(u)}{|u|^{3}}\,u,
\end{align*}
which concludes the proof of part (a).

We now turn to (b), and once again apply the decomposition \eqref{eq:k_decomp} of $S_{-}^{2}$:
\begin{align} \label{eq:b_a_lem_part(b)}
    \int_{S_{-}^{2}}&\,b_{a}\pfrac{\hat{k}\cdot u}{\eps}\hat{k}\otimes\hat{k}\,d\hat{k} 
    = \int_{0}^{1}b_{a}\pfrac{\lambda|u|}{\eps}\int_{S^{1}(u)}\paren{\sqrt{1-\lambda^{2}}\,\xi - \lambda\hat{u}} \otimes \paren{\sqrt{1-\lambda^{2}}\,\xi - \lambda\hat{u}}\,d\xi\,d\lambda.
\end{align}
We expand the tensor product using bilinearity and examine each term individually. First,
\begin{equation}
    \label{eq:b_a_lem_part(b)_1}
    \int_{S^{1}(u)}\sqrt{1-\lambda^{2}}\,\xi\otimes \sqrt{1-\lambda^{2}}\,\xi\,d\xi
    = (1-\lambda^{2})\int_{S^{1}(u)}\xi\otimes \xi\,d\xi,
\end{equation}
which by using the identity (stated and proved in Lemma \ref{lem:projection_matrix})  
\begin{equation} \label{eq:origin-proj-matrix-use}
  \int_{S^{1}(u)}\xi\otimes \xi\,d\xi =  
  \pi (I - \hat{u}\otimes \hat{u}), 
\end{equation}
becomes
\begin{equation*}
    \int_{S^{1}(u)}\sqrt{1-\lambda^{2}}\,\xi\otimes \sqrt{1-\lambda^{2}}\,\xi\,d\xi
    = \pi(1-\lambda^{2})(I - \hat{u}\otimes \hat{u}).
\end{equation*}
Meanwhile, the first cross-term in the expansion of \eqref{eq:b_a_lem_part(b)} vanishes:
\begin{align*}
    -\int_{S^{1}(u)}\sqrt{1-\lambda^{2}}\,\xi \otimes \lambda\hat{u}\,d\xi
    &= -\lambda\sqrt{1-\lambda^{2}}\paren{\int_{S^{1}(u)}\xi\,d\xi} \otimes \hat{u} = 0,
\end{align*}
and the other cross-term vanishes as well by the same reasoning. Finally,
\begin{align*}
    \int_{S^{1}(u)}\lambda\hat{u} \otimes \lambda\hat{u}\,d\xi
    = 2\pi\lambda^{2}\hat{u}\otimes \hat{u}.
\end{align*}
Consequently \eqref{eq:b_a_lem_part(b)} becomes
\begin{align}
    \int_{S_{-}^{2}}b_{a}\pfrac{\hat{k}\cdot u}{\eps}\hat{k}\otimes\hat{k}\,d\hat{k}
    &= \int_{0}^{1} b_{a}\pfrac{\lambda|u|}{\eps}
    \cdot \pi\brac{(1-\lambda^{2})(I-\hat{u}\otimes \hat{u}) + 2\lambda^{2}\hat{u}\otimes \hat{u}}\,d\lambda \nonumber \\
    &= \pi\int_{0}^{1} b_{a}\pfrac{\lambda|u|}{\eps}
    \brac{(I - \hat{u}\otimes \hat{u}) - \lambda^{2}(I - 3\hat{u}\otimes \hat{u})}\,d\lambda \label{eq:asymp_lemma_b_proof}
\end{align}
The first term in \eqref{eq:asymp_lemma_b_proof} becomes the main term of \eqref{eq:b_3} using the change of variables $\mu = \lambda|u|/\eps$:
\begin{align}
\pi\int_{0}^{1}b_{a}\pfrac{\lambda|u|}{\eps}(I - \hat{u}\otimes \hat{u})\,d\lambda
    &= \frac{\pi\eps}{|u|}\int_{0}^{|u|/\eps} b_{a}(\mu)(I - \hat{u}\otimes \hat{u})\,d\mu \nonumber \\
    &= \frac{8\pi^2\eps B_{a,\eps}(u)}{|u|}(I - \hat{u}\otimes \hat{u})
    \nonumber \\
    &=
    8\pi^2\eps B_{a,\eps}(u) K(u).  \label{eq:asymp_lemma_b_proof_change_var_line3}
\end{align}
For the second term of \eqref{eq:asymp_lemma_b_proof} we define
\begin{align}
    8\pi^2\eps R_\eps(u):=-
    \pi\int_0^1\! b_a\left( \frac{\lambda|u|}{\eps}\right)\lambda^2 (I - 3\hat u\otimes\hat u)\, d\lambda,
\end{align}
and apply the same change of variables as in \eqref{eq:asymp_lemma_b_proof_change_var_line3} to see that for any $\eta < 2$,
\begin{align*}
    \|R_{\eps}(u)\|
    &= \frac{\|I - 3\hat{u}\otimes \hat{u}\|}{8\pi\eps}\int_{0}^{|u|/\eps}\frac{\eps^{3}}{|u|^{3}}\mu^{a+2}\hat{V}(\mu)^{2}\,d\mu \\
    &=  \frac{\|I - 3\hat{u}\otimes \hat{u}\|}{8\pi\eps}\frac{\eps}{|u|}\int_{0}^{|u|/\eps}\pfrac{\eps\mu}{|u|}^{2}\mu^{a}\hat{V}(\mu)^{2}\,d\mu \\
    &\leq \frac{\|I - 3\hat{u}\otimes \hat{u}\|}{8\pi}\frac{1}{|u|}\int_{0}^{|u|/\eps}\pfrac{\eps\mu}{|u|}^{\eta}\mu^{a}\hat{V}(\mu)^{2}\,d\mu \\
    &\leq \frac{\|I - 3\hat{u}\otimes \hat{u}\|}{8\pi}M_{a+\eta}\,\frac{\eps^{\eta}}{|u|^{1+\eta}},
\end{align*}
where to obtain the third line we used the integration bound $\mu\le|u|/\eps$. Hence ($b$) is proved.

The same reasoning can be used to prove ($c$): if $\eta \leq 1$, then
\begin{align*}
    \int_{S_{-}^{2}}b_{a}\pfrac{\hat{k}\cdot u}{\eps}\,d\hat{k} 
    &= \int_{0}^{1}\int_{S^{1}(u)} b_{a}\pfrac{\lambda|u|}{\eps}\,d\xi\,d\lambda \\
    & = \frac{2\pi\eps}{|u|}\int_{0}^{|u|/\eps} \mu^{a}\hat{V}(\mu)^{2}\,d\mu \\
    &= 2\pi\pfrac{\eps}{|u|}^{\eta}\int_{0}^{|u|/\eps} \pfrac{\eps\mu}{|u|}^{1-\eta}\mu^{a-1+\eta}\hat{V}(\mu)^{2}\,d\mu \\
    &\leq 2\pi\pfrac{\eps}{|u|}^{\eta}\int_{0}^{|u|/\eps} \mu^{a-1+\eta}\hat{V}(\mu)^{2}\,d\mu \\
    &\leq 2\pi \pfrac{\eps}{|u|}^{\eta}M_{a-1+\eta}.
\end{align*}
\end{proof}

The other lemma of this section gives a similar bound for integrals involving
\begin{align}
    \label{eq:b_a_cross}
    b^{cr}_a\left(\frac{|u|}{\eps},\hat k\cdot \hat u\right) :=
    \frac{|\hat k\cdot u|^a}{\eps^a}\hat V\left(\frac{|\hat k\cdot u|}{\eps}\right)
			\hat V
			\left(\frac{\abs{u}}{\eps}\sqrt{1-(\hat k\cdot\hat u)^2}\right),
\end{align}
which arises in the analysis of the cross terms. This generalizes Lemma 4.1 of \cite{BenedettoPulvirenti07}.

\begin{lemma} \label{lem:b_a_cross} 
	For any $\eta\ge0$ and $a\ge 1$, we have the estimate
	\begin{align*}
		&\abs{\int_{S^2_-}\!b^{cr}_a\!\left(\frac{|u|}{\eps},\hat k\cdot \hat u\right)
			\, d\hat k }
		\lesssim_{a,\eta} M_{1+2\eta}
		\left( \frac{\eps}{|u|}\right)^{2+\eta-a}.
	\end{align*}
\end{lemma}
\begin{proof}
	We argue as in the proof of Lemma 4.1 of \cite{BenedettoPulvirenti07}. We define the variables $\gamma = |u|/\eps$ and $\lambda = |\hat u\cdot\hat k|$, and rewrite 
    \begin{align*}
        \abs{\int_{S^2_-}\!b^{cr}_a\!\left(\frac{|u|}{\eps},\hat k\cdot \hat u\right)
			\, d\hat k }
        &=
        2\pi\abs{
        \int_0^1\!  (\gamma\lambda)^a
        \hat V(\gamma\lambda)\hat V(\gamma\sqrt{1-\lambda^2})
        \, d\lambda
        }
        \\
        &\le
        2\pi
         \int_0^1\!  (\gamma\lambda)^a
        |\hat V(\gamma\lambda)||\hat V(\gamma\sqrt{1-\lambda^2})|
        \, d\lambda
        \\
        &=
        2\pi\int_{1/\sqrt2}^1\!\gamma^a\lambda(\lambda^{a-1} + \sqrt{1-\lambda^2}^{a-1})
        |\hat V(\gamma\lambda)||\hat V(\gamma\sqrt{1-\lambda^2})|
        \, d\lambda
        \\
        &\le
        C_a\gamma^a
        \int_{1/\sqrt2}^1\!\lambda
        |\hat V(\gamma\lambda)||\hat V(\gamma\sqrt{1-\lambda^2})|
        \, d\lambda
        \\
        &\le
        C_a\gamma^a
        \left(
            \int_{1/\sqrt2}^1\!\lambda|\hat V(\gamma\lambda)|^2\, d\lambda
        \right)^{1/2}
        \left(
            \int_0^{1/\sqrt2}\!\lambda|\hat V(\gamma\lambda)|^2\, d\lambda
        \right)^{1/2}
        \\
        &= C_a\gamma^{a-2}
        \left(
            \int_{\gamma/\sqrt2}^\gamma\!\mu|\hat V(\mu)|^2\, d\mu
        \right)^{1/2}
        \left(
            \int_0^{\gamma/\sqrt2}\!\mu|\hat V(\mu)|^2\, d\mu
        \right)^{1/2}
        \\
        &\le
        C_{a,\eta} \gamma^{a-2-\eta}
        \left(
            \int_{\gamma/\sqrt2}^\gamma\!\mu^{1+2\eta}|\hat V(\mu)|^2\, d\mu
        \right)^{1/2}
        \left(
            \int_0^{\gamma/\sqrt2}\!\mu|\hat V(\mu)|^2\, d\mu
        \right)^{1/2}
        \\
        &\le C_{a,\eta}\gamma^{a-2-\eta} \sqrt{I_1 I_{1+2\eta}},
    \end{align*} 
    where we define $I_a = \int_0^\infty\!\mu^a|\hat V(\mu)|\, d\mu$, and where we have freely used the change of variables $\lambda\mapsto\sqrt{1-\lambda^2}$, which preserves the measure $\lambda\,d\lambda$. Expanding $\gamma = |u|/\eps$ then proves the lemma.
\end{proof}

\section{Proof of Theorem \ref{mainthm}} \label{sect:proof_of_main}

In this section we prove the main result of this paper. 
We start by deriving a weak formulation 
\eqref{quantum_Landau} of the quantum Landau operator $Q_{qL}^{\alpha_0}(f)$. We then derive the weak formulation \eqref{eq:UU_weak_form} of the quantum Boltzmann operator $Q_{qB}^{\alpha,\eps}(f)$, which we will decompose into various terms that together converge to the weak form of $Q_{qL}^{\alpha_0}(f)$. Our method is inspired by the semiclassical limit result in \cite{BenedettoPulvirenti07}. However, we require additional Taylor expansions in $f$ to show convergence of the ternary term in the regime $\alpha(\eps)=O(1)$, and we consequently require additional derivative estimates on $f$ to close the limit.

\subsection{Weak Formulation of the Quantum Landau Operator} 
\label{subsec-qL-weak}

Before we formulate a result about the weak form of the operator 
$Q_{qL}^{\alpha_0}$ given by \eqref{quantum_Landau}, 
it will be useful to introduce the following notation. For any $f\in L^1_{\mathrm{\loc}}(\mathbb{R}^3)$ such that the convolution with the projection matrix $K$ (see \eqref{angular_projection}) is well-defined, we write
\begin{align}
    \label{eq:a_identity}
    a[f] := K*f
\end{align}
for the Landau diffusion matrix,
which satisfies the identities
\begin{align}
\label{eq:div_a_identity}
    \operatorname{div} a[f] = -\frac{2u}{|u|^3}*f
\end{align}
and
\begin{align*}
     -\operatorname{div}\operatorname{div}a[f] = 8\pi f.
\end{align*}
We can then rewrite the quantum Landau operator \eqref{quantum_Landau} in divergence form as
\begin{align}
    Q_{qL}^{\alpha_0}(f) 
    &=
    B\operatorname{div}_v\int_{\mathbb{R}^3}\! K(v-v_1)f(v_1)(1+\theta\alpha_0 f(v_1))\nabla_v f(v)\, d v_1
    \nonumber
    \\
    \nonumber
    &\qquad
    -
    B\operatorname{div}_v\int_{\mathbb{R}^3}\! \operatorname{div}K(v-v_1)f(v)(1+\theta\alpha_0 f(v))f(v_1)\, d v_1
    \\
    &= B\,\operatorname{div}\Big(a[f(1+\theta\alpha_0 f)]\nabla f - \operatorname{div}a[f]f(1+\theta\alpha_0 f)\Big),
    \label{quantum_Landau_convolution_div_form}
\end{align}
by integrating by parts to get the second integral and used standard properties of convolutions to move derivatives.

\begin{lemma}
\label{weak_formulation_quantum_Landau}
    Let $\psi\in W^{2,\infty}$ and $\langle v\rangle^N f\in W^{2,\infty}$ for some $N>3$. Then we have that 
    \begin{multline*}
        \int_{\R^3}\! 
        Q_{qL}^{\alpha_0}(f)(v)\psi(v)\,dv 
        \\
        = \int_{\R^3\times\R^3}\!
        \left[ - \frac{2B (v-v_1)}{|v-v_1|^3}\cdot(\nabla_v \psi -\nabla_{v_1}\psi_1) + B\Tr(K(v-v_1)D^2 \psi(v))\right] f\! f_1(1+\theta\alpha_0 f_1)\,dv_{1}\,dv,
    \end{multline*}
    for $Q_{qL}^{\alpha_0}$ as defined in \eqref{quantum_Landau}. Here we write $\psi_1$ for $\psi(v_1)$, $\psi$ for $\psi(v)$, and similarly for $f_1$ and $f$.
\end{lemma}
\begin{proof}
    We use \eqref{quantum_Landau_convolution_div_form} to rewrite
    \begin{align}
        \int_{\R^3}\! Q_{qL}^{\alpha_0}(f)(v)\psi(v)\,dv
        &= -B\int_{\R^3}\!\nabla_v \psi\cdot a[f(1+\theta\alpha_0 f)]\nabla f +
        B\int_{\mathbb{R}^3}\nabla_v \psi\cdot\operatorname{div} a[f]f(1+\theta\alpha_0 f)\,dv
        \label{eq:Qql_weak_integral}
        \\
        &=
        B\int_{\R^3}\!\Tr(D^2_v \psi(v) a[f(1+\theta\alpha_0 f)])f\,dv 
        \label{eq:Qql_weak_diffusion_integral-1}
        \\
        &\qquad\qquad+ B\int_{\R^3}\!\nabla_v \psi\cdot\big(\operatorname{div} a[f(1+\theta\alpha_0 f)]f+ \operatorname{div} a[f]f(1+\theta\alpha_0 f) \big)\,dv
        \label{eq:Qql_weak_diffusion_integral-2}
        \\
        \label{eq:Qql_weak_diffusion_integral-3}
        &=
        B\int_{\R^3\times\R^3}\!\left[\Tr K(v-v_1)D^2_v \psi\right]f\! f_1(1+\theta\alpha_0 f_1)
        \,dv_1\,dv
        \\
        \label{eq:Qql_weak_reaction_integral}
        &\qquad\qquad+
        B\int_{\R^3\times\R^3}\!\left[ -\frac{2(v-v_1)}{|v-v_1|^3}\cdot\nabla_v \psi\right] f\! f_1(2+\theta\alpha_0(f+f_1))\,dv_{1}\,dv,
        \end{align}
where to obtain \eqref{eq:Qql_weak_integral} we integrated by parts, to obtain \eqref{eq:Qql_weak_diffusion_integral-1} and 
\eqref{eq:Qql_weak_diffusion_integral-2}
we integrated by parts again in $v$ in the first integral, and to obtain final two lines we used the formulas \eqref{eq:a_identity} and \eqref{eq:div_a_identity} to rewrite the operators $a[\cdot]$ and $\operatorname{div} a[\cdot]$ as integrals in $v_1$. We can then rewrite the integral in \eqref{eq:Qql_weak_reaction_integral} as
\begin{align}
    B\int_{\R^3\times\R^3}\!&\left[ -\frac{2(v-v_1)}{|v-v_1|^3}\cdot\nabla_v \psi\right] f\! f_1(2+\theta\alpha_0(f+f_1))\,dv_{1}\,dv
    \nonumber \\
    &=
    \frac{B}{2}\int_{\R^3\times\R^3}\!\left[ -\frac{2(v-v_1)}{|v-v_1|^3}\cdot(\nabla_v \psi-\nabla_{v_1}\psi_1)\right] f\! f_1(2+\alpha_0 (f+f_1))\,dv
    _{1}\,dv
    \\
    &=
    B\int_{\R^3\times\R^3}\!\left[ -\frac{2(v-v_1)}{|v-v_1|^3}\cdot(\nabla_v \psi-\nabla_{v_1}\psi_1)\right] f\! f_1(1+\theta\alpha_0 f_1)\,dv
    _{1}\,dv,
    \label{eq:Qql_weak_reaction_integral_resymmetrized}
\end{align}
where we have first symmetrized in the variables $v$ and $v_1$, and then anti-symmetrized the $f\! f_1(2+\theta\alpha(f+f_1))$ term to match the (non-symmetric) expression $f\! f_1(1+\theta\alpha_0 f)$ in the integral in \eqref{eq:Qql_weak_diffusion_integral-3}.
Combining the term \eqref{eq:Qql_weak_diffusion_integral-3} with the term \eqref{eq:Qql_weak_reaction_integral} as rewritten in \eqref{eq:Qql_weak_reaction_integral_resymmetrized} then gives the desired formula
        \begin{multline*}
        \int_{\R^3}\! Q_{qL}^{\alpha_0}(f)(v)\psi(v)\,dv
        \\
        =B\int_{\R^3\times\R^3}\!\left[\Tr K(v-v_1)D^2_v \psi  - \frac{2(v-v_1)}{|v-v_1|^3}\cdot \left(\nabla_v \psi - \nabla_{v_1} \psi_1\right)\right]f\! f_1(1+\theta\alpha_0 f_1)
        \,dv_1\,dv.
    \end{multline*}
\end{proof}

\subsection{Weak Formulation of the Quantum Boltzmann Operator} \label{sect:expansion_of_uu}

In this section we justify our weak formulation \eqref{eq:UU_weak_form} of $Q_{qB}^{\alpha,\varepsilon}(f)$ by deriving it from our initial definition \eqref{eq:UU_formal_def}
\begin{multline*}
    Q_{qB}^{\alpha,\varepsilon}(f)(v)
    = \frac{1}{8\pi^{2}\eps^{2}}\int_{\R^{3}}\!\int_{\R^{3}}
    \brac{\hat{V}(k) + \theta\hat{V}\paren{k + \frac{u}{\eps}}}^{2}
    \delta(\eps |k|^{2} + k\cdot u) \\
    \{(1 + \theta\alpha f)(1 + \theta\alpha f_{1})f'\!f_{1}'
    - (1 + \theta\alpha f')(1 + \theta\alpha f'_{1})f\! f_{1}\}\,dk\,dv_{1}.
\end{multline*}
It will be useful to decompose $Q_{qB}^{\alpha,\varepsilon}(f)$ into the following form.

\begin{lemma}
\label{lemma:QUU_decomposition}
    If $\langle v\rangle ^N f\in L^\infty(\R^3)$ for some $N>3$, then
    \begin{align}
        \label{eq:QUU_decomposition}
        Q_{qB}^{\alpha,\varepsilon}(f)
    = Q_{B}^{\eps}(f) + Q_{T}^{\alpha,\eps}(f)
    + R_{B}^{\eps}(f) + R_{T}^{\alpha,\eps}(f),
    \end{align}
    where we define the terms
    \begin{align}
    \label{eq:Q_B_defn}
    Q_{B}^{\eps}(f)(v)
    &= \frac{1}{4\pi^{2}\eps^{4}}\int_{\R^{3}}\!\int_{S_{-}^{2}}
    |\hat{k}\cdot u|
    \hat{V}\pfrac{|\hat{k}\cdot u|}{\eps}^{2}
    (f'\!f_{1}' - f\! f_{1})\,d\hat{k}\,dv_{1} \\
    \label{Q_T_defn}
    Q_{T}^{\alpha,\eps}(f)(v)
    &= \frac{\theta\alpha}{4\pi^{2}\eps^{4}}\int_{\R^{3}}\!\int_{S_{-}^{2}} 
    |\hat{k}\cdot u|\hat{V}\pfrac{|\hat{k}\cdot u|}{\eps}^{2}
    \paren{f'\!f_{1}'(f+f_{1}) - f\! f_{1}(f'+f_{1}')}\,d\hat{k}\,dv_{1} \\
    \label{eq:R_B_defn}
    R_{B}^{\eps}(f)(v)
    &= \frac{\theta}{4\pi^{2}\eps^{4}}\int_{\R^{3}}\!\int_{S_{-}^{2}} 
    |\hat{k}\cdot u|\hat{V}\pfrac{|\hat{k}\cdot u|}{\eps}\hat{V}\Bigg(\frac{\sqrt{u^{2}-(\hat{k}\cdot u)^{2}}}{\eps}\Bigg)(f'\!f_{1}' - f\! f_{1})\,d\hat{k}\,dv_{1} \\
    \label{eq:R_T_defn}
    \begin{split}
    R_{T}^{\alpha,\eps}(f)(v)
    &= \frac{\alpha}{4\pi^{2}\eps^{4}}\int_{\R^{3}}\!\int_{S_{-}^{2}}
    |\hat{k}\cdot u|\hat{V}\pfrac{|\hat{k}\cdot u|}{\eps}
    \hat{V}\Bigg(\frac{\sqrt{u^{2}-(\hat{k}\cdot u)^{2}}}{\eps}\Bigg) \\
    &\qquad\qquad\qquad\qquad\qquad\qquad\qquad \cdot \paren{f'\!f_{1}'(f+f_{1}) - f\! f_{1}(f'+f_{1}')}\,d\hat{k}\,dv_{1},
    \end{split}
\end{align}
with $f' = f(v')$ and $f'_1 = f(v_1')$ now defined by the post-collisional velocities
\begin{align} \label{eq:post-collisional_velocities}
    v' &= v - (\hat{k}\cdot u)\hat{k} &
    v_1' &= v_1 + (\hat{k}\cdot u)\hat{k},
\end{align}
and the radial function $\hat{V}$ being viewed as a function on $[0,\infty)$.
\end{lemma}

\begin{remark} \label{remark:N>3}
The assumption that $\inner{v}^N f\in L^\infty$ for some $N>3$ that appears above and in Lemma \ref{lemma:weak_form} does not figure into either lemma's proof. Rather, it is a sufficient condition to ensure that the terms $Q_{B}^{\eps}(f)$, $Q_{T}^{\alpha,\eps}(f)$, $R_{B}^{\eps}(f)$, and $R_{T}^{\alpha,\eps}(f)$ are well-defined objects, specifically $L^1$ functions; this follows from the same reasoning as in the proof of the functional estimate Lemma \ref{lem:QUU_well-defined} (see Subsection \ref{subsect:proof_of_functional_estimate}). The functional estimate also requires a condition on $\hat{V}$, which we are assuming implicitly.
\end{remark}

\begin{proof}
For brevity, we will write
\begin{align}
\begin{split}
    \label{eq:defn_of_F}
    F(v,v_{1},k) &= (1 + \theta\alpha f)(1 + \theta\alpha f_{1})f'\!f_{1}'
    - (1 + \theta\alpha f')(1 + \theta\alpha f'_{1})f\! f_{1} \\
    &= (f'\!f_{1}' - f\! f_{1}) + \theta\alpha[f'\!f_{1}'(f + f_{1}) - f\! f_{1}(f'+f_{1}')],
\end{split}
\end{align}
where the dependence on $k$ comes from the terms $f' = f(v')$ and $f_{1}' = f(v_{1}')$ through the post-collisional velocities
\begin{align}
    \label{eq:post_collisional_velocities_k}
    v' &= v + \eps k & 
    v_{1}' &= v_{1} - \eps k
\end{align}
defined in \eqref{eq:postcollisional_velocities}. We start by expanding the potential part in $Q_{qB}^{\alpha,\varepsilon}$, to get
\begin{equation} \label{eq:expand_square}
\begin{aligned}
    Q_{qB}^{\alpha,\varepsilon}(f)(v) = 
    {}&\frac{1}{8\pi^{2}\eps^{2}}\int_{\R^{3}}\!\int_{\R^{3}}\hat{V}(k)^{2}
    \delta(\eps |k|^{2}+k\cdot u)F(v,v_{1},k)\,dk\,dv_{1} \\
    {}&+ \frac{1}{8\pi^{2}\eps^{2}}\int_{\R^{3}}\!\int_{\R^{3}}\hat{V}\paren{k + \frac{u}{\eps}}^{2}\delta(\eps |k|^{2} + k\cdot u)F(v,v_{1},k)\,dk\,dv_{1} \\
    {}&+\frac{\theta}{4\pi^{2}\eps^{2}}\int_{\R^{3}}\!\int_{\R^{3}}
    \hat{V}(k)\hat{V}\paren{k + \frac{u}{\eps}}
    \delta(\eps |k|^{2} + k\cdot u)F(v,v_{1},k)\,dk\,dv_{1}.
\end{aligned}
\end{equation}
In the second term on the right, we make the change of variables 
$$
k \mapsto k_{new}: =  -\paren{k + \frac{u}{\eps}}.
$$
This change of variables does not change the delta function, since
\begin{align*}
    \eps |k|^{2} + k\cdot u
    = \eps |k_{new}|^{2} + 2k_{new}\cdot u + \frac{u^{2}}{\eps}
    - k_{new}\cdot u - \frac{u^{2}}{\eps} 
    = \eps |k_{new}|^{2} + k_{new}\cdot u.
\end{align*}
Also, since 
\begin{align*}
    v':=v + \eps k
    &= v_{1} - \eps k_{new} &
    v_1':=v_{1} - \eps k
    &= v + \eps k_{new},
\end{align*}
this change of variables swaps the roles of $v'$ and $v_1'$, so $F(v,v_{1},k) = F(v,v_1,k_{new})$. Finally, the change of variables replaces $\hat{V}\paren{k + \frac{u}{\eps}}$ with $\hat{V}(-k_{new}) = \hat{V}(k_{new})$. We thus see that the second term of \eqref{eq:expand_square} is equal to the first, and so
\begin{align*}
    Q_{qB}^{\alpha,\varepsilon}(f)(v) = 
    {}&\frac{1}{4\pi^{2}\eps^{2}}\int_{\R^{3}}\!\int_{\R^{3}}\hat{V}(k)^{2}
    \delta(\eps |k|^{2}+k\cdot u)F(v,v_{1},k)\,dk\,dv_{1} \\
    {}&+\frac{\theta}{4\pi^{2}\eps^{2}}\int_{\R^{3}}\!\int_{\R^{3}}\hat{V}(k)\hat{V}\paren{k + \frac{u}{\eps}}
    \delta(\eps |k|^{2} + k\cdot u)F(v,v_{1},k)\,dk\,dv_{1}.
\end{align*}
The next step is to integrate the delta function using the following identity (Lemma \ref{lemma:delta}):
\begin{equation*}
    \int_{\R^{3}}\delta(\eps |k|^{2} + k\cdot u)g(k)\,dk
    = \frac{1}{\eps^{2}}\int_{S^{2}_{-}} |\hat{k}\cdot u|g\paren{-\frac{\hat{k}\cdot u}{\eps}\hat{k}}\,d\hat{k}.
\end{equation*}
Applying this and recasting the radial function $\hat{V}$ as a function on $[0,\infty)$ gives us
\begin{align*} 
    Q_{qB}^{\alpha,\varepsilon}(f)(v) &= 
    \frac{1}{4\pi^{2}\eps^{4}}\int_{\R^{3}}\!\int_{S_{-}^{2}}\! |\hat{k}\cdot u|\,\hat{V}\paren{\frac{|\hat{k}\cdot u|}{\eps}}^{2}F\left(v,v_{1},-\frac{\hat{k}\cdot u}{\eps}\hat{k}\right)\,d\hat{k}\,dv_{1} \\
    &\quad+\frac{\theta}{4\pi^{2}\eps^{4}}\int_{\R^{3}}\!\int_{S_{-}^{2}}\! |\hat{k}\cdot u|
    \hat{V}\paren{\frac{|\hat{k}\cdot u|}{\eps}}\hat{V}\Bigg(\frac{\sqrt{u^{2}-(\hat{k}\cdot u)^{2}}}{\eps}\Bigg)F\left(v,v_{1},-\frac{\hat{k}\cdot u}{\eps}\hat{k}\right)\,d\hat{k}\,dv_{1}.
\end{align*}
Substituting $k = -(\hat k\cdot u)\hat k/\varepsilon$ into the post-collisional velocities defined in \eqref{eq:post_collisional_velocities_k}
gives us 
\begin{align*}
    v' &= v - (\hat{k}\cdot u)\hat{k} &
    v_{1}' &= v_{1} + (\hat{k}\cdot u)\hat{k}
\end{align*}
as in the statement of the lemma, which concludes the proof of 
\eqref{eq:QUU_decomposition}.
\end{proof}

We now apply these terms to a test function, and use symmetrization to put the resulting expressions in forms which will prove useful when computing their limits.

\begin{lemma} \label{lemma:weak_form}
Let $\psi\in L^\infty(\R^3)$ and $ \langle v \rangle^N f\in L^\infty(\R^3)$ for some $N>3$. Then, for the operators $Q_{B}^{\eps}$, $Q_{T}^{\alpha,\eps}$, $R_{B}^{\eps}$, and $R_{T}^{\eps}$ defined in Lemma \ref{lemma:QUU_decomposition},
\begin{align}
    \inner{Q_{B}^{\eps}(f),\psi}
    &= \frac{1}{4\pi^{2}\eps^{4}}\int_{\R^{3}}\!\int_{\R^{3}}f\! f_{1}
    \int_{S_{-}^{2}}
    |\hat{k}\cdot u|
    \hat{V}\pfrac{|\hat{k}\cdot u|}{\eps}^{2}
    [\psi(v') - \psi(v)]\,d\hat{k}\,dv_{1}\,dv \label{eq:Q_MB} \\
    \inner{Q_{T}^{\alpha,\eps}(f),\psi}
    &= \frac{\theta\alpha}{4\pi^{2}\eps^{4}}\int_{\R^{3}}\!\int_{\R^{3}}f\! f_{1}\int_{S_{-}^{2}}
    |\hat{k}\cdot u|\hat{V}\pfrac{|\hat{k}\cdot u|}{\eps}^{2}
    [\psi(v')-\psi(v)](f'+f_{1}')\,d\hat{k}\,dv_{1}\,dv
    \label{eq:Q_T} \\
    \inner{R_{B}^{\eps}(f),\psi}
    &= \frac{\theta}{4\pi^{2}\eps^{4}}\int_{\R^{3}}\!\int_{\R^{3}}f\! f_{1}\int_{S_{-}^{2}}
    |\hat{k}\cdot u|\hat{V}\pfrac{|\hat{k}\cdot u|}{\eps} 
    \hat{V}\Bigg(\frac{\sqrt{u^{2}-(\hat{k}\cdot u)^{2}}} {\eps}\Bigg)
    [\psi(v') - \psi(v)]\,d\hat{k}\,dv_{1}\,dv
    \label{eq:R_B} \\ 
    \inner{R_{T}^{\alpha,\eps}(f),\psi}
    &= \frac{\alpha}{4\pi^{2}\eps^{4}}\int_{\R^{3}}\!\int_{\R^{3}}f\! f_{1}\int_{S_{-}^{2}}
    |\hat{k}\cdot u|\hat{V}\pfrac{|\hat{k}\cdot u|}{\eps}
    \hat{V}\Bigg(\frac{\sqrt{u^{2}-(\hat{k}\cdot u)^{2}}}{\eps}\Bigg) 
    [\psi(v') - \psi(v)](f'+f_{1}')\,d\hat{k}\,dv_{1}\,dv \label{eq:R_T}
\end{align}
with the post-collisional velocities $v'$ and $v_1'$ defined by 
\begin{align*}
 v' &= v - (\hat{k}\cdot u)\hat{k} &
    v_{1}' &= v_{1} + (\hat{k}\cdot u)\hat{k}
\end{align*}
as in \eqref{eq:post-collisional_velocities}.
\end{lemma}
\begin{remark}
For ease of notation, going forward we will write integration of a function $F\in L^1(\R^3)$ against a test function $\psi\in L^\infty(\R^3)$ as a pairing: $$\langle F,\psi\rangle := \int_{\R^3}\! F(v)\psi(v)\, dv.$$
As noted in Remark \ref{remark:N>3}, the terms $Q_{B}^{\eps}(f)$, $Q_{T}^{\alpha,\eps}(f)$, $R_{B}^{\eps}(f)$, and $R_{T}^{\alpha,\eps}(f)$ are $L^1$ when $\inner{v}^{N}f\in L^\infty$ for some $N>3$.
\end{remark}

\begin{proof}
We will treat the case $\inner{Q_{B}^{\eps}(f),\psi}$ explicitly; the other equations follow from the same reasoning. From the definition of $Q_{B}^{\eps}(f)$ in \eqref{eq:Q_B_defn}, our initial weak form expression is
\begin{equation} \label{eq:initial_weak_form_QB}
    \inner{Q_{B}^{\eps}(f),\psi}
    = \frac{1}{4\pi^{2}\eps^{4}}\int_{\R^{3}}\psi(v)\int_{\R^{3}}\!\int_{S_{-}^{2}}|\hat{k}\cdot u|\hat{V}\pfrac{|\hat{k}\cdot u|}{\eps}^{2}(f'\!f_{1}' - f\!f_{1})\,d\hat{k}\,dv_{1}\,dv.
\end{equation}
This can be written as a difference of two integrals, and we will use a change of variables to rewrite the first one,
\begin{equation} \label{eq:first_integral_of_QB}
    \int_{\R^{3}} \psi(v)\int_{\R^{3}}
    \int_{S_{-}^{2}}
    |\hat{k}\cdot u|
    \hat{V}\pfrac{|\hat{k}\cdot u|}{\eps}^{2}
    f'\!f_{1}'\,d\hat{k}\,dv_{1}\,dv.
\end{equation}
Specifically, we apply the pre-postcollisional change of variables
\begin{align*}
    v_{new} &= v' = v - (\hat{k}\cdot u)\hat{k} &
    v_{1,new} &= v_{1}' = v_{1} + (\hat{k}\cdot u)\hat{k}.
\end{align*}
Under this change of variables
\begin{gather*}
    u_{new}
    = v_{new} - v_{new}'
    = u - 2(\hat{k}\cdot u)\hat{k}, \\
    \hat{k}\cdot u_{new} = \hat{k}\cdot [u - 2(\hat{k}\cdot u)\hat{k}]
    = -\hat{k}\cdot u,
\end{gather*}
so we can rewrite \eqref{eq:first_integral_of_QB} as
\begin{align*}
    \int_{\R^{3}}\!\int_{\R^{3}}\!\int_{S_{-}^{2}}
    \psi(v')|\hat{k}\cdot u|
    \hat{V}\pfrac{|\hat{k}\cdot u|}{\eps}^{2}
    f\! f_{1}\,d\hat{k}\,dv_{1}\,dv.
\end{align*}
Substituting this back into \eqref{eq:initial_weak_form_QB} gives us
\begin{equation*}
    \inner{Q_{B}^{\eps}(f),\psi}
    = \frac{1}{4\pi^{2}\eps^{4}}\int_{\R^{3}}\!\int_{\R^{3}}f\! f_{1}
    \int_{S_{-}^{2}}
    |\hat{k}\cdot u|
    \hat{V}\pfrac{|\hat{k}\cdot u|}{\eps}^{2}
    [\psi(v') - \psi(v)]\,d\hat{k}\,dv_{1}\,dv,
\end{equation*}
as desired. The only difference when treating the other terms, namely $\inner{R_{B}^{\eps}(f),\psi}$ and $\inner{R_{T}^{\alpha,\eps}(f),\psi}$, is that one must consider how the change of variables affects the expression $\sqrt{u^{2} - (\hat{k}\cdot u)^{2}}$. However, since
\begin{equation*}
    u_{new}^{2}
    = u^{2} - 4(\hat{k}\cdot u)\hat{k}\cdot u + 4(\hat{k}\cdot u)^{2}
    = u^{2}
\end{equation*}
and $\hat{k}\cdot u_{new} = -\hat{k}\cdot u$, the expression remains unchanged. 
\end{proof}

In the subsequent sections, we will show that $Q_{B}^{\eps}$ combined with $Q_{T}^{\alpha,\eps}$ (the ``main terms'') converges to our desired kernel $Q_{qL}^{\alpha_0}$, while $R_{B}^{\eps}$ and $R_{T}^{\alpha,\eps}$ (the ``cross terms'') converge to 0.

\subsection{The Binary Main Term $Q_{B}^{\eps}$ } \label{sect:binary_main}

We will henceforth assume that $\psi\in W^{3,\infty}(\R^3)$, as in the statement of Theorem \ref{mainthm}. Starting from the weak form \eqref{eq:Q_MB}
\begin{equation*} 
    \inner{Q_{B}^{\eps}(f),\psi}
    = \frac{1}{4\pi^{2}\eps^{4}}\int_{\R^{3}}\!\int_{\R^{3}}f\! f_{1}
    \int_{S_{-}^{2}}
    |\hat{k}\cdot u|
    \hat{V}\pfrac{|\hat{k}\cdot u|}{\eps}^{2}
    [\psi(v') - \psi(v)]\,d\hat{k}\,dv_{1}\,dv,
\end{equation*}
we follow Benedetto-Pulvirenti \cite[p. 914]{BenedettoPulvirenti07} and rewrite $\psi(v') - \psi(v)$ using the Taylor expansion
\begin{equation}
\label{u'-u_Taylor}
    \psi(v') - \psi(v)
    = |\hat{k}\cdot u|\hat{k}\cdot \nabla_{v}\psi
    + \frac{1}{2}|\hat{k}\cdot u|^{2}\Tr(\hat{k}\otimes \hat{k}\,D^{2}\psi(v)) + |\hat{k}\cdot u|^{3}r_{3}(v,v'),
\end{equation}
where $\|r_3\|_{L^{\infty}(\R^{3}\times \R^{3})} \lesssim \|\psi\|_{W^{3,\infty}}$. (See Lemma \ref{lemma:taylor} for a proof.) Substituting \eqref{u'-u_Taylor} into \eqref{eq:Q_MB} results in the decomposition
\begin{equation} \label{eq:expansion_of_QB}
    \inner{Q_{B}^{\eps}(f),\psi}
    = \int_{\R^{3}}\!\int_{\R^{3}}f\! f_{1}(T_{2} + T_{3} + T_{4})\,dv_{1}\,dv,
\end{equation}
where
\begin{align}
    T_{2} &= \frac{1}{4\pi^{2}\eps^{2}}\int_{S_{-}^{2}}b_{2}\pfrac{|\hat{k}\cdot u|}{\eps}\hat{k}\cdot \nabla_{v}\psi(v)\,d\hat{k} ,
    \nonumber
    \\
    T_{3} &= \frac{1}{8\pi^{2}\eps}\int_{S_{-}^{2}}b_{3}\pfrac{|\hat{k}\cdot u|}{\eps}\Tr(\hat{k}\otimes \hat{k}\,D^{2}\psi(v))\,d\hat{k} ,
    \nonumber
    \\
    T_{4} &= \frac{1}{4\pi^{2}}\int_{S_{-}^{2}}b_{4}\pfrac{|\hat{k}\cdot u|}{\eps}r_{3}\,d\hat{k}.
    \label{eq:T4_defn}
\end{align}
(These were denoted in \cite{BenedettoPulvirenti07} as $T_{1}$, $T_{2}$, $T_{3}$.) Here we used the functions $b_a$ defined in \eqref{eq:def_of_b_a} and the identity
\begin{equation*}
    |\hat{k}\cdot u|^{a} \hat V\pfrac{|\hat{k}\cdot u|}{\eps}^{2}
    = \eps^{a}b_{a}\pfrac{|\hat{k}\cdot u|}{\eps}.
\end{equation*}

\subsection{The Ternary Main Term $Q_{T}^{\alpha,\eps}$} \label{sect:ternary_main}

Similarly to the decomposition of $Q_{B}^{\eps}$, we will decompose the ternary main term
\begin{equation} \label{eq:Q_T_before_taylor}
    \inner{Q_{T}^{\alpha,\eps}(f),\psi}
    = \frac{\theta\alpha}{4\pi^{2}\eps^{4}}\int_{\R^{3}}\!\int_{\R^{3}}f\! f_{1}\int_{S_{-}^{2}}
    |\hat{k}\cdot u|\hat{V}\pfrac{|\hat{k}\cdot u|}{\eps}^{2}
    [\psi(v')-\psi(v)](f'+f_{1}')\,d\hat{k}\,dv_{1}\,dv
\end{equation}
using the Taylor expansion in $\psi(v')-\psi(v)$ \eqref{u'-u_Taylor}.
However, in order to capture the ternary term in the limit, we also use a Taylor expansion of $f'+f_{1}'$:
\begin{align} \label{eq:f_Taylor}
    f' + f_{1}'
    = f + f_{1}
    + |\hat{k}\cdot u|\hat{k}\cdot (\nabla_{v}f - \nabla_{v_{1}}f_{1})
    + |\hat{k}\cdot u|^{2}\rho_{2}(v,v'),
\end{align}
where $\|\rho_2\|_{L^\infty(\R^3\times \R^3)} \lesssim \|f\|_{W^{2,\infty}}$ (see Lemma \ref{lemma:taylor} for a proof). Note that we have used that
\begin{equation*}
    v' - v = -(\hat{k}\cdot u)\hat{k}
    = |\hat{k}\cdot u|\hat{k}, 
\end{equation*}
which follows from the formula for $v'$ \eqref{eq:post-collisional_velocities} and the fact that $\hat{k}\in S_{-}^2$. We can now substitute \eqref{u'-u_Taylor} and \eqref{eq:f_Taylor} into the formula \eqref{eq:Q_T_before_taylor} for $\inner{Q_{T}^{\alpha,\eps}(f),\psi}$, resulting in
\begin{equation} \label{eq:expansion_of_QT}
    \inner{Q_{T}^{\alpha,\eps}(f),\psi} = \int_{\R^{3}}\!\int_{\R^{3}}f\! f_{1}(T_{2}' + T_{3a}' + T_{3b}' + T_{4}' + T_{5}' + T_{6}')\,dv_{1}\,dv.
\end{equation}
The terms $T_2'$, $T_{3a}'$, $T_{3b}'$ are
\begin{align}
    T_{2}'
    &= \frac{\theta\alpha}{4\pi^{2}\eps^{2}}\int_{S_{-}^{2}}
    b_{2}\pfrac{|\hat{k}\cdot u|}{\eps}
    (\hat{k}\cdot \nabla_{v}\psi)(f+f_{1})\,d\hat{k} ,
    \nonumber
    \\
    T_{3a}' &= \frac{\theta\alpha}{8\pi^{2}\eps}\int_{S_{-}^{2}}
    b_{3}\pfrac{|\hat{k}\cdot u|}{\eps}\Tr(\hat{k}\otimes \hat{k}\,D^{2}\psi(v))(f + f_{1})\,d\hat{k} ,
    \nonumber
    \\
    \label{eq:T_3b'}
    T_{3b}'
    &= \frac{\theta\alpha}{4\pi^{2}\eps}\int_{S_{-}^{2}}
    b_{3}\pfrac{|\hat{k}\cdot u|}{\eps}(\hat{k}\cdot \nabla_{v}\psi)\sparen{\hat{k}\cdot (\nabla_{v}f - \nabla_{v_{1}}f_{1})}\,d\hat{k},
\end{align}
and these are the terms which have non-trivial limits. The remaining terms are
\begin{align}
    \label{eq:T4'_defn}
    \begin{split}
    T_{4}'
    &= \frac{\theta\alpha}{4\pi^{2}}\int_{S_{-}^{2}}
    b_{4}\pfrac{|\hat{k}\cdot u|}{\eps}
    \bigg(r_{3}(f+f_{1}) + \rho_{2}(\hat{k}\cdot \nabla_{v}\psi) \\
    &\qquad\qquad\qquad\qquad\qquad\qquad+ \frac{1}{2}\Tr(\hat{k}\otimes \hat{k}\,D^{2}\psi(v))\hat{k}\cdot (\nabla_{v}f - \nabla_{v_{1}}f_{1})\bigg)\,d\hat{k} ,
    \end{split}
    \\
    \label{eq:T5'_defn}
    T_{5}'
    &= \frac{\theta\alpha \eps}{4\pi^{2}}\int_{S_{-}^{2}}
    b_{5}\pfrac{|\hat{k}\cdot u|}{\eps}
    \paren{\frac{1}{2}\,\rho_{2}\Tr(\hat{k}\otimes \hat{k}\,D^{2}\psi(v)) + r_{3}\hat{k}\cdot (\nabla_{v}f - \nabla_{v_{1}}f_{1})}\,d\hat{k} ,
    \\
    \label{eq:T6'_defn}
    T_{6}'
    &= \frac{\theta\alpha \eps^2}{4\pi^{2}}
    \int_{S_{-}^{2}}
    b_{6}\pfrac{|\hat{k}\cdot u|}{\eps}\rho_{2}r_{3}\,d\hat{k},
\end{align}
which will all converge to 0.

\subsection{Limits of the Main Terms} \label{sect:limits_of_main_terms}

Throughout this section we will encounter expressions of the following forms:
\begin{equation} \label{eq:integrals}
    \begin{aligned}
        I_{\lambda}^{(1)} &:= \int_{\R^3}\!\int_{\R^3}\frac{|f\!f_{1}|}{|u|^{\lambda}}\,dv_{1}\,dv \\
        I_{\lambda}^{(2)} &:= \int_{\R^3}\!\int_{\R^3}\frac{|f^2 f_{1}|}{|u|^{\lambda}}\,dv_{1}\,dv
        = \int_{\R^3}\!\int_{\R^3}\frac{|f\!f_{1}^2|}{|u|^{\lambda}}\,dv_{1}\,dv \\
        I_\lambda^{(3)} &:= \int_{\R^3}\!\int_{\R^3}\frac{|f\!f_{1}||\nabla_{v}f|}{|u|^{\lambda}}\,dv_{1}\,dv
        = \int_{\R^3}\!\int_{\R^3}\frac{|f\!f_{1}||\nabla_{v_{1}}f_{1}|}{|u|^{\lambda}}\,dv_{1}\,dv,
    \end{aligned}
\end{equation}
where $\lambda\in \R$. We can bound these integrals as follows:

\begin{lemma} \label{lemma:I_bounds}
If $\inner{v}^{N}f\in W^{1,\infty}(\R^3)$ for $N> 3$, then for all $\lambda\in (0,3)$,
\begin{align} \label{eq:I_bounds}
    I_{\lambda}^{(1)} &\lesssim_{N,\lambda} \|\inner{v}^{N}f\|_{L^{\infty}}^2 &
    I_{\lambda}^{(2)} &\lesssim_{N,\lambda} \|\inner{v}^{N}f\|_{L^{\infty}}^3 &
    I_{\lambda}^{(3)} &\lesssim_{N,\lambda} \|\inner{v}^{N}f\|_{W^{1,\infty}}^{3}.
\end{align}
\end{lemma}

\begin{proof}
A key observation is that if $\inner{v}^{N}f\in L^{\infty}(\R^3)$ for $N>3$, then $\inner{v}^{-N}\in L^p$ for all $1\leq p\leq\infty$, and so
\begin{equation*}
    \|f\|_{L^p} \leq \|\inner{v}^{N}f\|_{L^\infty}\|\inner{v}^{-N}\|_{L^p} 
    \lesssim_{N,p} \|\inner{v}^{N}f\|_{L^\infty},
    \qquad 1\leq p\leq \infty.
\end{equation*}
Therefore it suffices to get bounds in terms of $L^p$ norms of $f$ (or $\nabla f$, for $I_{\lambda}^{(3)}$). We start with $I_{\lambda}^{(1)}$. For any $\lambda \in (0,3)$, there exist $1<p,q<\infty$ (depending on $\lambda$) such that $I_{\lambda}^{(1)} \lesssim_\lambda \|f\|_{L^p}\|f\|_{L^q}$ by the bilinear form of the Hardy--Littlewood--Sobolev inequality (see Sobolev \cite[\S 4]{Sobolev}), which proves the desired bound on $I_{\lambda}^{(1)}$ by our prior observation. The bounds on $I_{\lambda}^{(2)}$ and $I_{\lambda}^{(3)}$ now follow from those on $I_{\lambda}^{(1)}$:
\begin{align*}
    I_{\lambda}^{(2)} &\leq \|f\|_{L^\infty}I_{\lambda}^{(1)}
    \lesssim_{N,\lambda} \|\inner{v}^{N}f\|_{L^\infty}^3 \\
    I_{\lambda}^{(3)}
    &\leq \|\nabla f\|_{L^\infty}I_{\lambda}^{(1)}
    \lesssim_{N,\lambda} \|\inner{v}^{N}\nabla f\|_{L^\infty}\|\inner{v}^{N}f\|_{L^\infty}^2
    \lesssim_{N} \|\inner{v}^{N}f\|_{W^{1,\infty}}^3.
\end{align*}
\end{proof}

\begin{remark}
This is not a sharp range of $\lambda$ for which such bounds hold (for example, the bound on $I_{\lambda}^{(1)}$ easily extends to $\lambda=0$), but this is all that we will need.
\end{remark}

\subsubsection{ The limit for $T_{2}+T_{2}'$} 

\begin{lemma} \label{lemma:limit_of_T2}
Let
\begin{equation*}
    \int_{\R^3}\!\int_{\R^3}f\!f_{1}(T_{2} + T_{2}')\,dv_{1}\,dv 
    = \frac{1}{4\pi^{2}\eps^{2}}\int_{\R^{3}}\!\int_{\R^{3}}\brac{\int_{S_{-}^{2}}b_{2}\pfrac{|\hat{k}\cdot u|}{\eps}\hat{k}\,d\hat{k}}\cdot (\nabla_{v}\psi)f\! f_{1}\sparen{1 + \theta\alpha(f+f_{1})}\,dv_{1}\,dv
\end{equation*}
and 
\begin{equation*}
    L_{2}(f) := 
    -B\int_{\R^{3}}\!\int_{\R^{3}}\paren{\frac{2u}{|u|^{3}}\cdot (\nabla_{v}\psi - \nabla_{v_{1}}\psi_{1})}f\! f_{1}\sparen{1 + \theta\alpha_{0}(f+f_{1})}\,dv_{1}\,dv.
\end{equation*}
Then for all $\eta \in (0,2)$, 
\begin{multline*}
    \abs{\int_{\R^3}\!\int_{\R^3}f\!f_{1}(T_{2} + T_{2}')\,dv_{1}\,dv
    - L_{2}(f)}  \\
    \lesssim_{N,\eta} \sparen{\eps^{\eta}M_{3+\eta}(1+\alpha_0)
        + M_3|\alpha-\alpha_{0}|}\|\psi\|_{W^{3,\infty}}\paren{\|\inner{v}^{N}f\|_{L^\infty}^2 + \|\inner{v}^{N}f\|_{L^\infty}^3}.
\end{multline*}
\end{lemma}

\begin{remark}
Note that in this lemma and for the rest of the paper, we will be taking
\begin{equation*}
    B := \int_{0}^{\infty}\mu^3 \hat{V}(\mu)^2\,d\mu < M_3.
\end{equation*}
\end{remark}

\begin{proof}
We start by using the asymptotic lemma (Lemma \ref{lemma:b_a}(a)) and symmetrization to write
\begin{multline} \label{eq:T2-T2_rewrite}
    \frac{1}{4\pi^{2}\eps^{2}}\int_{\R^{3}}\!\int_{\R^{3}}\brac{\int_{S_{-}^{2}}b_{2}\pfrac{|\hat{k}\cdot u|}{\eps}\hat{k}\,d\hat{k}}\cdot (\nabla_{v}\psi)f\! f_{1}\sparen{1 + \theta\alpha(f+f_{1})}\,dv_{1}\,dv \\
    \begin{aligned} 
        &= -4\int_{\R^{3}}\!\int_{\R^{3}}\paren{\frac{B_{3,\eps}(u)u}{|u|^{3}}\cdot \nabla_{v}\psi}f\! f_{1}\sparen{1 + \theta\alpha(f+f_{1})}\,dv_{1}\,dv \\
        &= -2\int_{\R^{3}}\!\int_{\R^{3}}B_{3,\eps}(u)\paren{\frac{u}{|u|^{3}}\cdot (\nabla_{v}\psi - \nabla_{v_{1}}\psi_{1})}f\! f_{1}\sparen{1 + \theta\alpha(f+f_{1})}\,dv_{1}\,dv.
    \end{aligned}
\end{multline}
Therefore
\begin{multline*}
    \int_{\R^3}\!\int_{\R^3}f\!f_{1}(T_{2} + T_{2}')\,dv_{1}\,dv
    - L_{2}(f) \\
    \begin{aligned}
        &= -\int_{\R^{3}}\!\int_{\R^{3}}\sparen{B_{3,\eps}(u) - B}\paren{\frac{2u}{|u|^{3}}\cdot (\nabla_{v}\psi - \nabla_{v_{1}}\psi_{1})}f\! f_{1}\sparen{1 + \theta\alpha_0(f+f_{1})}\,dv_{1}\,dv \\
        &\qquad {}- \theta(\alpha-\alpha_{0})\int_{\R^3}\!\int_{\R^{3}}B_{3,\eps}(u)\paren{\frac{2u}{|u|^{3}}\cdot (\nabla_{v}\psi - \nabla_{v_{1}}\psi_{1})}f\! f_{1}(f+f_1)\,dv_{1}\,dv.
    \end{aligned}
\end{multline*}
To bound this, note that 
\begin{align} \label{eq:bound_on_grad_psi_diff}
    \abs{\frac{u}{|u|^{3}}\cdot (\nabla_{v}\psi - \nabla_{v_{1}}\psi_{1})}
    \leq \frac{|\nabla_{v}\psi - \nabla_{v_1}\psi_1|}{|u|^2}
    \lesssim \frac{\|\psi\|_{W^{2,\infty}}}{|u|},
\end{align}
and
\begin{align} \label{eq:bound_on_B_diff}
    |B_{3,\eps}(u) - B|
    = \frac{1}{8 \pi}\int_{|u|/\eps}^{\infty}\mu^3 \hat{V}(\mu)^{2}\,d\mu
    \leq \frac{1}{8\pi}\frac{\eps^\eta}{|u|^\eta}\int_{|u|/\eps}^{\infty}\mu^{3+\eta}\hat{V}(\mu)^{2}\,d\mu 
    \leq \frac{1}{8\pi}\frac{\eps^{\eta}}{|u|^{\eta}}M_{3+\eta}
\end{align}
for all $\eta \geq 0$. It now follows from the integral bounds \eqref{eq:I_bounds} that if $\eta\in (0,2)$ (so that $\eta,1+\eta\in (0,3)$), then:
\begin{multline*}
    \abs{\int_{\R^3}\!\int_{\R^3}f\!f_{1}(T_{2} + T_{2}')\,dv_{1}\,dv
    - L_{2}(f)} \\
    \begin{aligned}
        &\lesssim \eps^{\eta}M_{3+\eta}\|\psi\|_{W^{3,\infty}}\paren{I_{1+\eta}^{(1)} + I_{\eta}^{(1)} + \alpha_{0}I_{1+\eta}^{(2)} + \alpha_{0}I_{\eta}^{(2)}}
        + B|\alpha - \alpha_{0}|\|\psi\|_{W^{3,\infty}}\paren{I_{1+\eta}^{(2)} + I_{\eta}^{(2)}} \\
        &\lesssim_{N,\eta} \sparen{\eps^{\eta}M_{3+\eta}(1+\alpha_0)
        + M_3|\alpha-\alpha_{0}|}\|\psi\|_{W^{3,\infty}}\paren{\|\inner{v}^{N}f\|_{L^\infty}^2 + \|\inner{v}^{N}f\|_{L^\infty}^3}.
    \end{aligned}
\end{multline*}
\end{proof}

\subsubsection{ The limit for $T_3 + T_{3a}'$}

\begin{lemma} \label{lemma:limit_of_T3a}
Let
\begin{multline*}
    \int_{\R^{3}}\!\int_{\R^{3}}f\! f_{1}(T_{3} + T_{3a}')\,dv_{1}\,dv \\
    = \frac{1}{8\pi^{2}\eps}\int_{\R^{3}}\!\int_{\R^{3}}\Tr\brac{\paren{\int_{S_{-}^{2}}b_{3}\pfrac{|\hat{k}\cdot u|}{\eps}\hat{k}\otimes \hat{k}\,d\hat{k}}D^{2}\psi}
    f\! f_{1}\sparen{1 
    +\theta\alpha(f+f_{1})}\,dv_{1}\,dv
\end{multline*}
and 
\begin{equation*}
    L_{3a}(f) := B\int_{\R^{3}}\!\int_{\R^{3}} \Tr[K(u)D^{2}\psi]f\! f_{1}\sparen{1 + \theta\alpha_0(f+f_{1})}\,dv_{1}\,dv.
\end{equation*}
Then for all $\eta \in (0,2)$,
\begin{multline*}
    \abs{\int_{\R^3}\!\int_{\R^3}f\!f_{1}(T_{3} + T_{3a'})\,dv_{1}\,dv
    - L_{3a}(f)} \\
    \lesssim_{N,\eta} \sparen{\eps^{\eta}M_{3+\eta}(1+\alpha_0 + \alpha) + M_{3}|\alpha-\alpha_0|}\|\psi\|_{W^{2,\infty}}
        \paren{\|\inner{v}^{N}f\|_{L^\infty}^2 + \|\inner{v}^{N}f\|_{L^\infty}^3}.
\end{multline*}
\end{lemma}
\begin{proof}
We start by using the asymptotic lemma (Lemma \ref{lemma:b_a}(b)) to write
\begin{multline*}
    \frac{1}{8\pi^{2}\eps}\int_{\R^{3}}\!\int_{\R^{3}}\Tr\brac{\paren{\int_{S_{-}^{2}}b_{3}\pfrac{|\hat{k}\cdot u|}{\eps}\hat{k}\otimes \hat{k}\,d\hat{k}}D^{2}\psi}
    f\! f_{1}\sparen{1 
    +\theta\alpha(f+f_{1})}\,dv_{1}\,dv \\
    = \int_{\R^{3}}\!\int_{\R^{3}} \Tr\brac{(B_{3,\eps}(u)K(u) + R_{\eps})D^{2}\psi}f\! f_{1}\sparen{1 + \theta\alpha(f+f_{1})}\,dv_{1}\,dv,
\end{multline*}
with
\begin{equation} \label{eq:bound_on_R_eps_in_T3a}
    \|R_{\eps}(u)\| \lesssim M_{3+\eta}\frac{\eps^{\eta}}{|u|^{1+\eta}},
    \qquad \eta < 2.
\end{equation}
Therefore 
\begin{align*}
    \int_{\R^3}\!\int_{\R^3}f\!f_{1}(T_{3} + T_{3a'})\,dv_{1}\,dv
    - L_{3a}(f) 
        &= \int_{\R^{3}}\!\int_{\R^{3}} \Tr\brac{(B_{3,\eps}(u)-B)K(u)D^{2}\psi}f\! f_{1}\sparen{1 + \theta\alpha_0(f+f_{1})}\,dv_{1}\,dv \\
        &\qquad + \theta(\alpha-\alpha_0)\int_{\R^{3}}\!\int_{\R^{3}} \Tr\brac{B_{3,\eps}(u)K(u)D^{2}\psi}f\! f_{1}(f+f_{1})\,dv_{1}\,dv \\
        &\qquad + \int_{\R^{3}}\!\int_{\R^{3}} \Tr\brac{R_\eps D^{2}\psi}f\! f_{1}\sparen{1 + \theta\alpha(f+f_{1})}\,dv_{1}\,dv.
\end{align*}
Using \eqref{eq:bound_on_R_eps_in_T3a}, the bound \eqref{eq:bound_on_B_diff} on $|B_{3,\eps}(u) - B|$ from the proof of Lemma \ref{lemma:limit_of_T2}, the integral bounds \ref{eq:I_bounds}, and 
\begin{equation} \label{eq:bound_on_K}
    \|K(u)\| \leq \frac{\|I\| + \|\hat{u}\otimes \hat{u}\|}{|u|}
    \lesssim \frac{1 + |\hat{u}|^{2}}{|u|}
    \lesssim \frac{1}{|u|},
\end{equation}
we see that
\begin{multline*}
    \abs{\int_{\R^3}\!\int_{\R^3}f\!f_{1}(T_{3} + T_{3a'})\,dv_{1}\,dv
    - L_{3a}(f)} \\
    \begin{aligned}
        &\lesssim \eps^{\eta}M_{3+\eta}\|\psi\|_{W^{2,\infty}}\paren{I_{1+\eta}^{(1)} + \alpha_{0}I_{1+\eta}^{(2)}}
        + B|\alpha-\alpha_0|\|\psi\|_{W^{2,\infty}}I_{1}^{(2)}
        + \eps^{\eta}M_{3+\eta}\|\psi\|_{W^{2,\infty}}\paren{I_{1+\eta}^{(1)} + \alpha I_{1+\eta}^{(2)}} \\
        &\lesssim_{N,\eta} \sparen{\eps^{\eta}M_{3+\eta}(1+\alpha_0 + \alpha) + M_{3}|\alpha-\alpha_0|}\|\psi\|_{W^{2,\infty}}
        \paren{\|\inner{v}^{N}f\|_{L^\infty}^2 + \|\inner{v}^{N}f\|_{L^\infty}^3}
    \end{aligned}
\end{multline*}
for all $\eta \in (0,2)$.
\end{proof}

\subsubsection{The limit for $T_{3b}'$} \label{subsect:T3b'}

\begin{lemma} \label{lemma:limit_of_T3b}
Let 
\begin{equation*}
    \int_{\R^3}\!\int_{\R^3}f\!f_{1}T_{3b}'\,dv_{1}\,dv 
    = \frac{\theta\alpha}{4\pi^{2}\eps}\int_{\R^3}\!\int_{\R^3}
    \brac{\int_{S_{-}^{2}}
    b_{3}\pfrac{|\hat{k}\cdot u|}{\eps}(\hat{k}\cdot \nabla_{v}\psi)\sparen{\hat{k}\cdot (\nabla_{v}f - \nabla_{v_{1}}f_{1})}\,d\hat{k}}f\!f_{1}\,dv_{1}\,dv
\end{equation*}
and 
\begin{equation*}
    L_{3b}(f) := 
    -\theta\alpha_{0} B\int_{\R^{3}}\!\int_{\R^{3}}\brac{-\frac{2u}{|u|^{3}}\cdot (\nabla_{v}\psi - \nabla_{v_{1}}\psi_1) + \Tr(K(u) D^{2}\psi)}f^2 f_{1}\,dv_{1}\,dv.
\end{equation*}
Then for all $\eta \in (0,2)$, 
\begin{equation*}
    \abs{\int_{\R^3}\!\int_{\R^3}f\!f_{1}T_{3b}'\,dv_{1}\,dv - L_{3b}(f)}  \lesssim_{N,\eta}
    \sparen{\eps^{\eta}M_{3+\eta}(\alpha_0 +\alpha) + M_{3}|\alpha-\alpha_{0}|}\|\psi\|_{W^{1,\infty}}\|\inner{v}^{N}f\|_{W^{1,\infty}}^{3}.
\end{equation*}
\end{lemma}
\begin{proof}
We will first write $L_{3b}(f)$ in a form more comparable to the $T_{3b}'$ integral. Using the formula \eqref{eq:div_a_identity} for $\div_{v}K(u)$, integration by parts, and symmetrization, we have
\begin{align*}
    L_{3b}(f) 
    &= 
    -\theta\alpha_{0} B\int_{\R^{3}}\!\int_{\R^{3}}\brac{-\frac{2u}{|u|^{3}}\cdot (\nabla_{v}\psi - \nabla_{v_{1}}\psi_1) + \Tr(K(u) D^{2}\psi)}f^2 f_1\,dv_{1}\,dv 
    \\
    &=
    -\theta\alpha_{0} B\int_{\R^{3}}\!\int_{\R^{3}}[(\div_{v} K(u))\cdot (\nabla_{v}\psi - \nabla_{v_{1}}\psi_{1}) +  \Tr(K(u)D^{2}\psi)]f^2 f_1\,dv_{1}\,dv 
    \\
    &=\theta\alpha_{0} B\int_{\R^{3}}\!\int_{\R^{3}}\Tr\brac{K(u)\paren{(\nabla_{v}\psi - \nabla_{v_{1}}\psi_{1}) \otimes \nabla_{v}(f^2 f_1)}}\,dv_{1}\,dv \\
    &= 2\theta\alpha_{0} B\int_{\R^{3}}\!\int_{\R^{3}}\Tr\brac{K(u)\paren{(\nabla_{v}\psi - \nabla_{v_{1}}\psi_{1}) \otimes \nabla_{v}f}}f\! f_1\,dv_{1}\,dv \\
    &= 2\theta\alpha_{0} B \int_{\R^{3}}\!\int_{\R^{3}}\Tr\brac{K(u)\paren{\nabla_{v}\psi \otimes (\nabla_{v}f - \nabla_{v_{1}}f_{1}}}f\! f_1\,dv_{1}\,dv.
\end{align*}
Likewise, using the identity
\begin{align*}
    (\xi\cdot x)(\xi\cdot y)
    = \Tr((\xi\otimes \xi)(x\otimes y))
    \qquad \xi,x,y\in \R^3
\end{align*}
and Lemma \ref{lemma:b_a}(b), we get
\begin{align*} 
    \int_{\R^3}\!\int_{\R^3}f\!f_{1}T_{3b}'\,dv_{1}\,dv 
    &= \frac{\theta\alpha}{4\pi^{2}\eps}\int_{\R^3}\!\int_{\R^3}
    \int_{S_{-}^{2}}
    \brac{b_{3}\pfrac{|\hat{k}\cdot u|}{\eps}(\hat{k}\cdot \nabla_{v}\psi)\sparen{\hat{k}\cdot (\nabla_{v}f - \nabla_{v_{1}}f_{1})}\,d\hat{k}}f\!f_{1}\,dv_{1}\,dv \\
    &= 2\theta\alpha\int_{\R^3}\!\int_{\R^3}\Tr\brac{\paren{B_{3,\eps}(u)K(u) + R_{\eps}(u)}\paren{\nabla_{v}\psi \otimes (\nabla_{v}f - \nabla_{v_{1}}f_{1})}}f\!f_{1}\,dv_{1}\,dv,
\end{align*}
where 
\begin{equation} \label{eq:bound_on_R_eps_in_T3b}
    \|R_{\eps}(u)\| \lesssim M_{3+\eta}\frac{\eps^{\eta}}{|u|^{1+\eta}},
    \qquad \eta < 2.
\end{equation}
Therefore
\begin{align*}
    \int_{\R^3}\!\int_{\R^3}f\!f_{1}T_{3b}'\,dv_{1}\,dv  - L_{3b}(f)
    &= 
    2\theta\alpha_{0}\int_{\R^3}\!\int_{\R^3}\Tr\brac{(B_{3,\eps}(u) - B)K(u)\paren{\nabla_{v}\psi \otimes (\nabla_{v}f - \nabla_{v_{1}}f_{1})}}f\!f_{1}\,dv_{1}\,dv \\
    &\qquad + 2\theta(\alpha-\alpha_{0})\int_{\R^3}\!\int_{\R^3}\Tr\brac{B_{3,\eps}(u)K(u)\paren{\nabla_{v}\psi \otimes (\nabla_{v}f - \nabla_{v_{1}}f_{1})}}f\!f_{1}\,dv_{1}\,dv \\
    &\qquad + 2\theta\alpha\int_{\R^3}\!\int_{\R^3}\Tr\brac{R_{\eps}\paren{\nabla_{v}\psi \otimes (\nabla_{v}f - \nabla_{v_{1}}f_{1})}}f\!f_{1}\,dv_{1}\,dv.
\end{align*}
Using \eqref{eq:bound_on_R_eps_in_T3b}, the bound \eqref{eq:bound_on_B_diff} on $|B_{3,\eps}(u) - B|$ from the proof of \ref{lemma:limit_of_T2}, the bound \eqref{eq:bound_on_K} on $K(u)$, and the integral bounds \eqref{eq:I_bounds}, we see that
\begin{align*}
    \abs{\int_{\R^3}\!\int_{\R^3}f\!f_{1}T_{3b}'\,dv_{1}\,dv  - L_{3b}(f)}
    &\lesssim
    \eps^{\eta}\alpha_0 M_{3+\eta}\|\psi\|_{W^{1,\infty}}I_{1+\eta}^{(3)}
    + B|\alpha-\alpha_{0}|\|\psi\|_{W^{1,\infty}}I_{1}^{(3)}
    + \alpha \eps^{\eta}M_{3+\eta}\|\psi\|_{W^{1,\infty}}I_{1+\eta}^{(3)} \\
    &\lesssim_{N,\eta} \sparen{\eps^{\eta}M_{3+\eta}(\alpha_0 +\alpha) + M_{3}|\alpha-\alpha_{0}|}\|\psi\|_{W^{1,\infty}}\|\inner{v}^{N}f\|_{W^{1,\infty}}^{3}.
\end{align*}
for all $\eta \in (0,2)$.
\end{proof}

\subsubsection{The remaining terms}  \label{subsect:remaining_terms}

\begin{lemma} \label{lemma:limit_of_T4_T5_T6}
Let 
\begin{align*}
    &\int_{\R^3}\!\int_{\R^3}f\!f_{1}(T_{4} + T_{4}')\,dv_{1}\,dv 
        =  \frac{1}{4\pi^{2}}\int_{\R^3}\!\int_{\R^3}
        \bigg[\int_{S_{-}^{2}} b_{4}\pfrac{|\hat{k}\cdot u|}{\eps}
        \bigg(r_{3}\sparen{1 + \theta\alpha(f+f_{1})} + \rho_{2}(\hat{k}\cdot \nabla_{v}\psi) \\
        &\qquad\qquad\qquad\qquad\qquad\qquad\qquad\qquad + \frac{1}{2}\Tr(\hat{k}\otimes \hat{k}\,D^{2}\psi(v))\hat{k}\cdot (\nabla_{v}f - \nabla_{v_{1}}f_{1})\bigg)\,d\hat{k}\bigg]f\!f_{1}\,dv_{1}\,dv \\
    &\int_{\R^3}\!\int_{\R^3}f\!f_{1}T_{5}'\,dv_{1}\,dv 
        = \frac{\theta\alpha \eps}{4\pi^{2}}
        \int_{\R^3}\!\int_{\R^3}\bigg[\int_{S_{-}^{2}}
        b_{5}\pfrac{|\hat{k}\cdot u|}{\eps}
        \bigg(\frac{1}{2}\,\rho_{2}\Tr(\hat{k}\otimes   \hat{k}\,D^{2}\psi(v))  \\
        &\qquad\qquad\qquad\qquad\qquad\qquad\qquad\qquad + r_{3}\hat{k}\cdot (\nabla_{v}f - \nabla_{v_{1}}f_{1})\bigg)\,d\hat{k}\bigg]f\!f_{1}\,dv_{1}\,dv \\
    &\int_{\R^3}\!\int_{\R^3}f\!f_{1}T_{6}'\,dv_{1}\,dv 
        = \frac{\theta\alpha \eps^2}{4\pi^{2}}
        \int_{\R^3}\!\int_{\R^3}\brac{\int_{S_{-}^{2}}
        b_{6}\pfrac{|\hat{k}\cdot u|}{\eps}\rho_{2}r_{3}\,d\hat{k}}f\!f_{1}\,dv_{1}\,dv.
\end{align*}
Then for all $\eta\in (0,1]$,
\begin{equation*}
    \abs{\int_{\R^3}\!\int_{\R^3}f\!f_{1}(T_{4} + T_{4}' + T_{5}' + T_{6}')\,dv_{1}\,dv} \lesssim_{N,\eta}
    \eps^{\eta}(M_{3+\eta} + \eps M_{4+\eta} + \eps^2 M_{5+\eta})(1+\alpha)\|\psi\|_{W^{3,\infty}}\|\inner{v}^{N}f\|_{W^{2,\infty}}^3.
\end{equation*}
\end{lemma}
\begin{proof}
We first apply the triangle inequality to each integral and apply the bounds $\|\rho_{2}\|_{L^\infty}\lesssim \|f\|_{W^{2,\infty}}$ and $\|r_{3}\|_{L^\infty} \lesssim \|\psi\|_{W^{3,\infty}}$ (proved in Lemma \ref{lemma:taylor}) to get
\begin{multline*}
    \abs{\int_{\R^3}\!\int_{\R^3}f\!f_{1}(T_{4} + T_{4}' + T_{5}' + T_{6}')\,dv_{1}\,dv} \\
    \begin{aligned}
    &\lesssim \int_{\R^3}\!\int_{\R^3}\bigg[\int_{S_{-}^2}
    b_{4}\pfrac{|\hat{k}\cdot u|}{\eps}
    \paren{\|\psi\|_{W^{3,\infty}}(1 + \alpha|f+f_{1}|) + \|f\|_{W^{2,\infty}}\|\psi\|_{W^{1,\infty}} + \|\psi\|_{W^{2,\infty}}\|f\|_{W^{1,\infty}}} \\
    &\qquad\qquad\qquad + b_{5}\pfrac{|\hat{k}\cdot u|}{\eps}\alpha\eps \paren{\|f\|_{W^{2,\infty}}\|\psi\|_{W^{2,\infty}} + \|\psi\|_{W^{3,\infty}}\|f\|_{W^{1,\infty}}} \\
    &\qquad\qquad\qquad + b_{6}\pfrac{|\hat{k}\cdot u|}{\eps}\alpha \eps^2 \|f\|_{W^{2,\infty}}\|\psi\|_{W^{3,\infty}}\bigg] 
    |f\!f_{1}|\,dv_{1}\,dv.
    \end{aligned}
\end{multline*}
With this, the only dependence on $\hat{k}$ in the integrand is in the $b_{a}$ functions. We can therefore apply Lemma \ref{lemma:b_a}(c) to each $b_a$ to get
\begin{multline*}
    \abs{\int_{\R^3}\!\int_{\R^3}f\!f_{1}(T_{4} + T_{4}' + T_{5}' + T_{6}')\,dv_{1}\,dv} \\
    \begin{aligned}
    &\lesssim
        \eps^{\eta}M_{3+\eta}\paren{\|\psi\|_{W^{3,\infty}}\paren{1 + \alpha \|f\|_{L^{\infty}}}I_{\eta}^{(1)}
        + \|f\|_{W^{2,\infty}}\|\psi\|_{W^{2,\infty}}I_{\eta}^{(1)}} \\
    &\qquad +  \eps^{\eta+1}\alpha M_{4+\eta}\|f\|_{W^{2,\infty}}\|\psi\|_{W^{3,\infty}}I_{\eta}^{(1)}
    + \eps^{\eta+2}\alpha M_{5+\eta}\|f\|_{W^{2,\infty}}\|\psi\|_{W^{3,\infty}}I_{\eta}^{(1)}
    \end{aligned}
\end{multline*}
for any $\eta \leq 1$. Using the integral bounds \eqref{eq:I_bounds}, which impose the bounds $\eta \in (0,3)$, we conclude 
\begin{align*}
    &\abs{\int_{\R^3}\!\int_{\R^3}f\!f_{1}(T_{4} + T_{4}' + T_{5}' + T_{6}')\,dv_{1}\,dv} \\
    &\qquad\qquad\qquad\qquad\lesssim_{N,\eta} 
    \eps^{\eta}(M_{3+\eta} + \eps M_{4+\eta} + \eps^2 M_{5+\eta})(1+\alpha)\|\psi\|_{W^{3,\infty}}\|f\|_{W^{2,\infty}}\|\inner{v}^{N}f\|_{L^\infty}^2 \\
    &\qquad\qquad\qquad\qquad\lesssim
    \eps^{\eta}(M_{3+\eta} + \eps M_{4+\eta} + \eps^2 M_{5+\eta})(1+\alpha)\|\psi\|_{W^{3,\infty}}\|\inner{v}^{N}f\|_{W^{2,\infty}}^3
\end{align*}
for all $\eta\in (0,1]$.
\end{proof}

\subsection{Limits of the Cross Terms} \label{sect:cross_terms}

It remains to prove convergence of the cross terms $\langle R_B^\eps(f),\psi\rangle$ and $\langle R_T^{\alpha,\eps}(f),\psi\rangle$. We start with the ternary cross term as defined in \eqref{eq:R_T}:
\begin{equation*}
    \langle R_T^{\alpha,\eps}(f),\psi\rangle
    = \frac{\alpha}{4\pi^{2}\eps^{4}}\int_{\R^{3}}\!\int_{\R^{3}}f\! f_{1}\int_{S_{-}^{2}}
    |\hat{k}\cdot u|\hat{V}\pfrac{|\hat{k}\cdot u|}{\eps}
    \hat{V}\Bigg(\frac{\sqrt{u^{2}-(\hat{k}\cdot u)^{2}}}{\eps}\Bigg) 
    [\psi(v') - \psi(v)](f'+f_{1}')\,d\hat{k}\,dv_{1}\,dv.
\end{equation*}
Our first step is to take Taylor expansions of $\psi(v') - \psi(v)$ and $f'+f_1'$ like in \eqref{u'-u_Taylor} and \eqref{eq:f_Taylor}, though to lower order. Specifically, we write
\begin{multline*}
	\langle R_T^{\alpha,\eps}(f),\psi\rangle
	=\frac{\alpha}{4\pi^2\eps^4}
	\int_{\R^{3}}\!\int_{\R^{3}}f\! f_{1}\int_{S_{-}^{2}}
	|\hat{k}\cdot u|\hat{V}\pfrac{|\hat{k}\cdot u|}{\eps}
	\hat{V}\Bigg(\frac{\sqrt{u^{2}-(\hat{k}\cdot u)^{2}}}{\eps}\Bigg) 
	\\
	\cdot
	\left(|\hat k\cdot u|\hat{k}\cdot (\nabla_v \psi(v) - \nabla_{v_1}\psi(v_1)) + |\hat{k}\cdot u|^2 \rho_2(v,v')\right)
	(f + f_1 + |\hat{k}\cdot u|\varrho_1(v,v'))\,d\hat{k}\,dv_{1}\,dv,
\end{multline*}
where
\begin{align*}
	\|\rho_{2}\|_{L^{\infty}} &\lesssim  \|\psi\|_{W^{2,\infty}},
	&
	\|\varrho_1\|_{L^\infty} &\lesssim  \|f \|_{W^{1,\infty}},
\end{align*}
(see Lemma \ref{lemma:taylor} for a proof) and therefore we can use Lemma \ref{lem:b_a_cross} to estimate
\begin{align}
    &\abs{\langle R_T^{\alpha,\eps}(f),\psi\rangle} \nonumber
    \\
    &\le
    \frac{\alpha}{4\pi^2}
    \int_{\mathbb R^3}\!\int_{\mathbb R^3}\! |f\! f_1|
    \int_{S^2_-}
    \|\psi\|_{W^{2,\infty}} \nonumber
    \\
    &\qquad\cdot
    \brac{
        \left(\frac{1}{\eps^2}b^{cr}_{2}
        \left(\frac{|u|}{\eps},\hat k\cdot u\right)\hat k + 
        \frac1\eps b^{cr}_{3}\left(\frac{|u|}{\eps},\hat k\cdot u\right)\right)|f+f_1|
        +
        b^{cr}_{4}\left(\frac{|u|}{\eps},\hat k\cdot u\right)\|f\|_{W^{1,\infty}}
    }\, d\hat k\, d v_1\, dv \nonumber
    \\
    &\lesssim \alpha \eps^{\eta'-2}M_{1+2\eta'}\|\psi\|_{W^{2,\infty}}
    \paren{I_{\eta'}^{(2)} + I_{\eta'-1}^{(2)} + I_{\eta'-2}^{(1)}\|f\|_{W^{1,\infty}}}
\end{align}
for any $\eta'\geq 0$, where the cross coefficients $b_a^{cr}$ were defined in \eqref{eq:b_a_cross}. We can then write $\eta = \eta'-2$ and apply the integral bounds \eqref{eq:integrals} to get
\begin{align}\label{eq:limit_of_RT}
    \abs{\langle R_T^{\alpha,\eps}(f),\psi\rangle}
    &\lesssim_{N,\eta} 
    \alpha\eps^\eta M_{5+2\eta}\|\psi\|_{W^{2,\infty}}\|\paren{\|\inner{v}^{N}f\|_{L^{\infty}}^3 + \|\inner{v}^{N}f\|_{L^{\infty}}^2 \|f\|_{W^{1,\infty}}} \nonumber \\
    &\lesssim_{N} \alpha\eps^\eta M_{5+2\eta}\|\psi\|_{W^{2,\infty}}\|\|\inner{v}^{N}f\|_{W^{1,\infty}}^3.
\end{align}
for any $\eta\in (0,1)$ (as we need $\eta,\eta+1,\eta+2\in (0,3)$ to apply \eqref{eq:integrals}). We proceed similarly for $R_{B}^{\eps}(f)$ except for not having the term $f'+f_{1}'$, and we see that
\begin{equation} \label{eq:limit_of_RB}
    \abs{\langle R_B^{\eps}(f),\psi\rangle} 
    \lesssim_{N,\eta} 
    \alpha\eps^\eta M_{5+2\eta} 
    \|\psi\|_{W^{2,\infty}}
     \|\inner{v}^{N}f\|_{W^{1,\infty}}^2.
\end{equation}

\subsection{Conclusion} \label{subsect:conclusion}

We can now finish the proof of Theorem \ref{mainthm}. Using the decomposition of the quantum Boltzmann operator in Lemmas \ref{lemma:QUU_decomposition}, the problem is reduced to showing that, for any test function $\psi$,
\begin{equation} \label{eq:overall_limit}
    \inner{Q_{B}^{\eps}(f)
    + Q_{T}^{\alpha,\eps}(f)
    + R_{B}^{\eps}(f)
    + R_{T}^{\eps}(f),\psi}
    \xto{\eps\to 0} \inner{Q_{qL}^{\alpha_{0}},\psi},
\end{equation}
where these terms are defined in Lemma \ref{lemma:weak_form}. As described in \eqref{eq:expansion_of_QB} and \eqref{eq:expansion_of_QT}, we can further decompose the two main terms as 
\begin{equation} \label{eq:QB_QT_summary}
    \inner{Q_{B}^{\eps}(f) + Q_{T}^{\alpha,\eps}(f),\psi}
    = \int_{\R^{3}}\!\int_{\R^{3}}f\! f_{1}\brac{(T_{2} + T_{2}') + (T_{3} + T_{3a}') + T_{3b}' + (T_{4} + T_{4}') + T_{5}' + T_{6}'}\,dv_{1}\,dv.
\end{equation}
Lemmas \ref{lemma:limit_of_T2}, \ref{lemma:limit_of_T3a}, \ref{lemma:limit_of_T3b}, and \ref{lemma:limit_of_T4_T5_T6} imply that
\begin{align*}
    &\int_{\R^{3}}\!\int_{\R^{3}}f\! f_{1}(T_{2}+T_{2}')\,dv_{1}\,dv
    \xto{\eps\to 0} B\int_{\R^{3}}\!\int_{\R^{3}}\paren{-\frac{2u}{|u|^{3}}\cdot (\nabla_{v}\psi -\nabla_{v_{1}}\psi_{1})}f\! f_{1}\sparen{1 + \theta\alpha_{0}(f+f_{1})}\,dv_{1}\,dv \\
    &\int_{\R^{3}}\!\int_{\R^{3}}f\! f_{1}(T_{3}+T_{3a}')\,dv_{1}\,dv
    \xto{\eps\to 0} B\int_{\R^{3}}\!\int_{\R^{3}} \Tr[K(u)D^{2}\psi]f\! f_{1}\sparen{1 + \theta\alpha_0(f+f_{1})}\,dv_{1}\,dv \\
    &\int_{\R^{3}}\!\int_{\R^{3}}f\! f_{1}T_{3b}'\,dv_{1}\,dv
    \xto{\eps\to 0} -\theta\alpha_{0} B\int_{\R^{3}}\!\int_{\R^{3}}\brac{-\frac{2u}{|u|^{3}}\cdot (\nabla_{v}\psi - \nabla_{v_{1}}\psi_1) + \Tr(K(u) D^{2}\psi)}f^2 f_{1}\,dv_{1}\,dv \\
    &\int_{\R^{3}}\!\int_{\R^{3}}f\! f_{1}(T_{4}+T_{4}'+T_{5}'+T_{6}')
    \xto{\eps\to 0} 0.
\end{align*}
Therefore
\begin{equation*}
    \inner{Q_{B}^{\eps}(f)+Q_{T}^{\alpha,\eps}(f),\psi} \xto{\eps\to 0} 
    \int_{\R^3\times\R^3}\!
        \left[ - \frac{2Bu}{|u|^3}\cdot(\nabla_v \psi -\nabla_{v_1}\psi_1) + B\Tr(K(u)D^2 \psi)\right] f\! f_1(1+\theta\alpha_0 f_1)\,dv_{1}\,dv,
\end{equation*}
which, as shown in Lemma \ref{weak_formulation_quantum_Landau}, is precisely $\inner{Q_{qL}^{\alpha_{0}}(f),\psi}$. Meanwhile, by \eqref{eq:limit_of_RB} and \eqref{eq:limit_of_RT} from Subsection \ref{sect:cross_terms},
\begin{equation*}
    \inner{R_{B}^{\eps}(f)+R_{T}^{\eps}(f),\psi} 
    \xto{\eps\to 0} 0,
\end{equation*}
which proves \eqref{eq:overall_limit}. Moreover, the convergence rates given by the aforementioned lemmas for the main terms are
\begin{align*}
    &\abs{\int_{\R^{3}}\!\int_{\R^{3}}f\! f_{1}(T_{2}+T_{2}')\,dv_{1}\,dv - L_{2}(f)}  \\
    &\qquad\qquad\qquad \lesssim_{N,\eta} \sparen{\eps^{\eta}M_{3+\eta}(1+\alpha_0)
        + M_3|\alpha-\alpha_{0}|}\|\psi\|_{W^{3,\infty}}\paren{\|\inner{v}^{N}f\|_{L^\infty}^2 + \|\inner{v}^{N}f\|_{L^\infty}^3} \\
    &\abs{\int_{\R^{3}}\!\int_{\R^{3}}f\! f_{1}(T_{3}+T_{3a}')\,dv_{1}\,dv - L_{3a}(f)} \\
    &\qquad\qquad\qquad \lesssim_{N,\eta} \sparen{\eps^{\eta}M_{3+\eta}(1+\alpha_0 + \alpha) + M_{3}|\alpha-\alpha_0|}\|\psi\|_{W^{2,\infty}}
        \paren{\|\inner{v}^{N}f\|_{L^\infty}^2 + \|\inner{v}^{N}f\|_{L^\infty}^3} \\
    &\abs{\int_{\R^{3}}\!\int_{\R^{3}}f\! f_{1}T_{3b}'\,dv_{1}\,dv - L_{3b}(f)}  \\
    &\qquad\qquad\qquad \lesssim_{N,\eta} \sparen{\eps^{\eta}M_{3+\eta}(\alpha_0 +\alpha) + M_{3}|\alpha-\alpha_{0}|}\|\psi\|_{W^{1,\infty}}\|\inner{v}^{N}f\|_{W^{1,\infty}}^{3} \\
    &\abs{\int_{\R^{3}}\!\int_{\R^{3}}f\! f_{1}(T_{4}+T_{4}'+T_{5}'+T_{6}')\,dv_{1}\,dv}
    \lesssim_{N,\eta} \eps^{\eta}(M_{3+\eta} + \eps M_{4+\eta} + \eps^2 M_{5+\eta})(1+\alpha)\|\psi\|_{W^{3,\infty}}\|\inner{v}^{N}f\|_{W^{2,\infty}}^3
\end{align*}
for any $N>3$ and $\eta\in (0,1]$ and the rates for the cross terms are
\begin{equation*}
    |\inner{R_{B}^{\eps}(f)+R_{T}^{\alpha,\eps}(f),\psi}|
    \lesssim 
    \alpha\eps^\eta M_{5+2\eta} 
    \|\psi\|_{W^{2,\infty}}
     \paren{\|\inner{v}^{N}f\|_{W^{1,\infty}}^2 + \|\inner{v}^{N}f\|_{W^{1,\infty}}^3}
\end{equation*}
for any $N>3$ and $\eta\in (0,1)$ by \eqref{eq:limit_of_RB} and \eqref{eq:limit_of_RT}. Since $M_{a}\leq M_{5+2\eta}$ for $a\in [1,5+2\eta]$, we conclude that if $\eta\in (0,1)$, $N>3$, and $\eps\in (0,1]$, then the overall convergence rate of the limit in \eqref{eq:overall_limit} is dominated by
\begin{equation*}
    (\eps^\eta + |\alpha-\alpha_0|)
    M_{5+2\eta}
    (1 + \alpha_0 + \alpha)
    \|\psi\|_{W^{3,\infty}}
    \paren{\|\inner{v}^{N}f\|_{W^{1,\infty}}^2 + \|\inner{v}^{N}f\|_{W^{2,\infty}}^3}.
\end{equation*}
Assuming all of the moments and norms are finite, this converges to 0 as $\eps\to 0$, completing the proof of the theorem.

\appendix

\section{Projection Matrix} \label{subsect:projection_matrix}

In this section we prove a basic identity that is used in the proof of Lemma \ref{lemma:b_a} in Subsection \ref{subsect:asymptotic_lemmas}. 

\begin{lemma}
    \label{lem:projection_matrix}
    Let $\hat u\in\mathbb{R}^3$ be such that $|\hat u|=1$, and define $S^1(\hat u) = \{\xi\in S^2\,:\,\xi\cdot\hat u=0\}$. Then we have the identity
\begin{equation} \label{eq:origin-proj-matrix}
  \int_{S^{1}(\hat u)}\xi\otimes \xi\,d\xi =  
  \pi (I - \hat{u}\otimes \hat{u}).
\end{equation}
\end{lemma}
\begin{proof}
If $A$ is an orthogonal matrix such that $Ae_3 = \hat{u}$ (where $e_3 = (0,0,1)$), then we can change variables $\zeta = A^{-1}\xi$ such that $\zeta\in S^1(e_3)$,
since $A^{-1}\xi\cdot e_3 = \xi\cdot \hat u=0$. We can then compute
\begin{equation}
    \label{eq:S1_xi_xi_integral}
    \int_{S^{1}(\hat u)}\xi\otimes \xi\,d\xi
    = \int_{S^{1}(e_3)}A\zeta \otimes A\zeta\,d\zeta
    =
    A\left(\int_{S^1(e_3)}\zeta\otimes\zeta\, d\zeta\right) A^{-1}.
\end{equation}
However, in polar coordinates, we can isolate the integral in \eqref{eq:S1_xi_xi_integral}. By computing 
\begin{align*}
    \int_{S^1(e_3)}\zeta\otimes\zeta\, d\zeta
    &=
    \int_0^{2\pi}\!
    \begin{pmatrix}
        \cos^2\vartheta & \cos\vartheta\sin\vartheta & 0\\
        \cos\vartheta\sin\vartheta & \sin^2\vartheta & 0\\
        0 & 0 & 0
    \end{pmatrix}d\vartheta
    =
    \pi\begin{pmatrix}1&0&0\\0&1&0\\0&0&0\end{pmatrix},
\end{align*}
we can therefore rewrite \eqref{eq:S1_xi_xi_integral} as 
\begin{align*}
    \pi A(I - e_3\otimes e_3)A^{-1}
    = \pi(I - (A e_3)\otimes (A e_3))
    =\pi( I - \hat u\otimes\hat u).
\end{align*}
\end{proof}

\section{Delta Function Identity} \label{subsect:delta_identity}

We here prove an identity which is used to expand the initial definition of the quantum Boltzmann operator, which is also found in \cite[(2.11)]{BenedettoPulvirenti07}. 

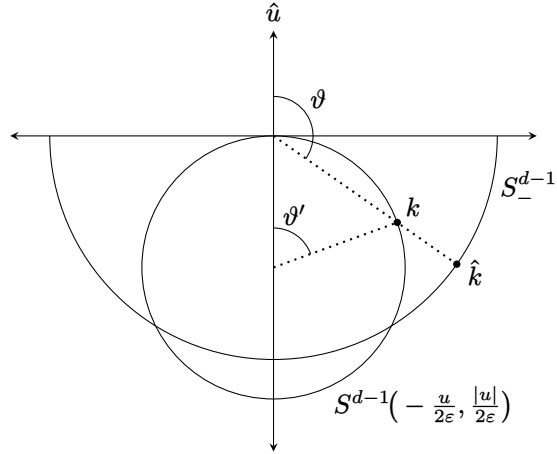
\begin{center}
\begin{tikzpicture}[scale=3.5, >=stealth]
	
	\draw[<->] (-1,0) -- (1,0);
	\draw[<->] (0,-1.2) -- (0,0.4) node[above] {$\hat u$};
	
	\draw (0,-0.5) circle (0.5);
	
	\coordinate (O) at (0,0);           
	\coordinate (C) at (0,-0.5);          
	\coordinate (P) at (0.47,-0.329);  
	\coordinate (B) at (0.696, -0.4875); 
	
	\draw[dotted,thick] (C) -- (P);
	\draw[dotted,thick] (O) -- (B);
	
	\draw(0,0) ++(90:0.15) arc[start angle=90, end angle=-35, radius=0.15];
	\node at (0.17,0.14) {$\vartheta$};
	
	\draw(0,-0.5) ++(90:0.15) arc[start angle=90, end angle=20, radius=0.15];
	\node at (0.08,-0.3) {$\vartheta'$};
	
	\draw (-0.85,0) arc[start angle=180, end angle=360, radius=0.85];
	
	\fill (P) circle (0.014);
	\node at (0.53,-0.27) {$k$};
	
	\fill (B) circle (0.014);
	\node at (0.766, -0.5175) {$\hat k$}; 
	
	\node at (0.97, -0.2) {$S^{d-1}_-$};
	
	\node at (0.57,-1.02) {$S^{d-1}\!\big(-\frac{u}{2\eps},\frac{|u|}{2\eps} \big)$};
	
\end{tikzpicture}
\end{center}

\begin{lemma} \label{lemma:delta}
Let $g$ be a bounded Borel-measurable function on $\R^{d}$ for any dimension $d\ge 2$, and let $u\in \R^{d}\sm \{0\}$. Then
\begin{equation} \label{eq:delta}
    \int_{\R^{d}}\delta(\eps |k|^{2} + k\cdot u)g(k)\,dk
    = \frac{1}{\eps^{d-1}}\int_{S^{d-1}_{-}} |\hat{k}\cdot u|^{d-2}g\paren{-\frac{\hat{k}\cdot u}{\eps}\hat{k}}\,d\hat{k},
\end{equation}
where
\begin{align*}
	S^{d-1}_-:= \{\hat k\in S^{d-1} :  \hat k\cdot u\le 0\}.
\end{align*}
\end{lemma}

We will provide two proofs of this Lemma, with the first using the specific geometry of the problem to take a more straightforward path, while the second is a change of variables argument with differential forms, which we hope may prove useful in other problems.

Before looking at either proof, we unpack the meaning of the left side of \eqref{eq:delta}, which both explains what this identity means for general $g$, and forms the starting point of both proofs. 

Since the derivative of the function $k\mapsto  \eps |k|^{2} + k\cdot u$ from $\R^3\to \R$ is nonzero on the zero set
\begin{equation} \label{eq:delta_zero_set}
    \wt{M} := \{k\in \R^d : \eps |k|^2 + k\cdot u = 0\},
\end{equation}
$\wt{M}$ is a hypersurface in $\R^d$, and $\delta(\eps|k|^2 + k\cdot u)$ is a distribution acting on test functions by
\begin{equation} \label{eq:delta_identity_unpacked}
    \int_{\R^d}\delta(\eps|k|^2 + k\cdot u)g(k)\,dk
    = \int_{\wt{M}}\frac{g(k)}{|2\eps k + u|}\,dV(k),
    \qquad g\in C_{c}^{\infty}(\R^d)
\end{equation}
where $dV$ denotes the surface measure on $\wt{M}$ (see \cite[Theorem 6.1.5]{Hormander}). Using this formula to rewrite the \eqref{eq:delta} results in an identity that makes sense even if $g$ is bounded Borel-measurable on $\R^d$ since, as we will see in the first proof, the manifold $\wt{M}$ is compact.

It will also be convenient for both proofs to simplify the right side of \eqref{eq:delta_identity_unpacked} by noting that if $k\in \wt{M}$, then
\begin{align*}
	\abs{2\eps k + u}^2
	= 4\eps^2 	\abs{k+\frac{u}{2\eps}}^2
    = 4\eps^2 \paren{\frac{\eps |k|^2 + k\cdot u}{\eps} + \frac{|u|^2}{4\eps^2}}
	= 4\eps^2 \frac{\abs{u}^2}{4\eps^2} =|u|^2.
\end{align*}
Therefore, the goal of both proofs becomes showing that
\begin{equation} \label{eq:delta_id_goal}
    \frac{1}{|u|}\int_{\wt{M}}g(k)\,dV(k)
    = \frac{1}{\eps^{d-1}}\int_{S^{d-1}_{-}} |\hat{k}\cdot u|^{d-2}g\paren{-\frac{\hat{k}\cdot u}{\eps}\hat{k}}\,d\hat{k}.
\end{equation}

\subsection{First Proof of Lemma \ref{lemma:delta}}

\begin{proof}
This proof starts by noting that the zero set $\wt{M}$ \eqref{eq:delta_zero_set} is actually a sphere:
\begin{equation*}
    \wt{M}
    = \{k : \eps |k|^2 + k\cdot u = 0\}
    = \braces{k : \eps\paren{\abs{k + \frac{u}{2\eps}}^2 - \frac{|u|^2}{4\eps^2}} = 0}
    = S^{d-1}\paren{-\frac{u}{2\eps},\frac{|u|}{2\eps}}.
\end{equation*}
Therefore, if we first assume that $u$ is of the form $|u|e_{d}$, then we can use the change of variables
\begin{align*}
	h = \frac{2\eps}{|u|}k+e_d
\end{align*}
to rewrite the left side of \eqref{eq:delta_id_goal} as
\begin{align}
	\label{eq:delta_identity_eq2}
	\frac{|u|^{d-2}}{(2\eps)^{d-1}}
	\int_{S^{d-1}}\! g\left( \frac{|u|}{2\eps}(h-e_d)\right)\,d\mathcal H^{d-1}(h),
\end{align}
where $\mathcal{H}^{d-1}$ denotes the $(d-1)$-dimensional Hausdorff measure on $S^{d-1}$. Re-expressing $h\in\R^d$ in hyperspherical coordinates as 
$$h = (r\sin(\vartheta')\omega,r\cos(\vartheta')),\qquad
(r,\omega,\vartheta')\in [0,\infty)\times S^{d-2} \times [0,\pi],$$
we can rewrite the integral \eqref{eq:delta_identity_eq2} as
\begin{align}
	&\frac{|u|^{d-2}}{(2\eps)^{d-1}}
	\int_0^\pi\!
	\int_{S^{d-2}}\! g\left( \frac{|u|}{2\eps}\sin(\vartheta')\omega,\frac{|u|}{2\eps}\left(\cos\vartheta'-1\right)\right)\sin^{d-2}(\vartheta')\,d\mathcal H^{d-2}(\omega)\,d\vartheta'
		\nonumber
	\\
	&
	=
	\frac{|u|^{d-2}}{(2\eps)^{d-1}}
	\int_{\frac\pi2}^\pi\!
	\int_{S^{d-2}}\!
		g\left( \frac{|u|}{\eps}\abs{\cos\vartheta}\sin(\vartheta)\omega,\frac{|u|}{\eps}\abs{\cos\vartheta}\cos\vartheta\right)
		2^{d-1}\abs{\cos\vartheta}^{d-2}\sin^{d-2}(\vartheta)\,d\mathcal H^{d-2}(\omega)\, 
		\,d\vartheta
	\tag{$\vartheta' = 2\vartheta-\pi$}
	\\
	&=
	\frac{1}{\eps^{d-1}}
			\int_{\frac\pi2}^\pi\!
		\int_{S^{d-2}}\!
		|u\cdot\hat k|^{d-2}
		g\left( \frac{|u|}{\eps}\abs{\cos\vartheta}\sin(\vartheta)\omega,\frac{|u|}{\eps}\abs{\cos\vartheta}\cos\vartheta\right)
		\sin^{d-2}(\vartheta)\,d\mathcal H^{d-2}(\omega)\, 
		\,d\vartheta
	\nonumber
	\\
	&=
	\frac{1}{\eps^{d-1}}
	\int_{\frac\pi2}^\pi\!
	\int_{S^{d-2}}\!
	|u\cdot\hat k|^{d-2}
	g\left( 
	-\frac{\hat k\cdot u}{\eps}
	\sin(\vartheta)\omega,	-\frac{\hat k\cdot u}{\eps}\cos\vartheta\right)
	\sin^{d-2}(\vartheta)\,d\mathcal H^{d-2}(\omega)\, 
	\,d\vartheta
	\nonumber
	\\
	&=
	\frac{1}{\eps^{d-1}}\int_{S^{d-1}_{-}} |\hat{k}\cdot u|^{d-2}g\paren{-\frac{\hat{k}\cdot u}{\eps}\hat{k}}\,d\hat{k},
	\nonumber
\end{align}
where we have used that
\begin{align*}
	\frac{|u|}{\eps}\abs{\cos\vartheta} = -\frac{u\cdot\hat k}{\eps},\qquad\frac\pi2\le\vartheta\le\pi,
\end{align*}
and thus the identity \eqref{eq:delta} is proved whenever $u = |u| e_d$. Rotational symmetry of \eqref{eq:delta} then proves the identity for all $u\in\R^d\setminus\{0\}$, which concludes the proof.
\end{proof}

\subsection{Second Proof of Lemma \ref{lemma:delta}}

This proof relies on the following well-known change of variables formula, for which we include a proof for completeness.

\begin{lemma} \label{lemma:change_of_var}
Let $M_{1},M_{2}$ be oriented $m$-dimensional submanifolds of $\R^{d}$, $g:M_{2}\to \R$ a measurable function, $F$ a smooth function from an open neighborhood of $M_{1}$ to $\R^{d}$ that maps $M_{1}$ diffeomorphically onto $M_{2}$, and let $U$ be an open subset of $M_1$ with smooth parametrization $\psi:V\subset \R^m \to U$. Then
\begin{equation*}
    \int_{F(U)}g\, dV_{2}
    = \int_{U}g(F(y))\paren{\frac{\det[(D\psi_{\psi^{-1}(y)})^{t}(DF_{y})^{t}(DF_{y})(D\psi_{\psi^{-1}(y)})]}{\det[(D\psi_{\psi^{-1}(y)})^{t}(D\psi_{\psi^{-1}(y)})]}}^{1/2}\,dV_{1}(y)
\end{equation*}
holds whenever either integral is absolutely convergent, where $dV_{i}$ denotes the volume form on $M_{i}$.
\end{lemma}
\begin{proof}
We first note that, given a measurable function $h:M_{1}\to \R$, we have
\begin{align*}
    \int_{U}h\,dV_{1}
    = \int_{\psi(V)}h\,dV_{1}
    = \int_{V}\psi^{*}(h\,dV_{1})
    = \int_{V}(h\circ \psi)\,\psi^{*}dV_{1}.
\end{align*}
The function $\psi^{-1}=(y^{1},\ldots, y^{k})$ gives a set of local coordinates on $U$, which lets us write
\begin{equation*}
    dV_{1} = \det(g_{ij})^{1/2}\,dy^{1}\wedge \cdots \wedge dy^{k},
\end{equation*}
where $g$ is the metric on $V_{1}$. Thus, computing the pullback $\psi^{*}dV_{1}$ boils down to computing $g_{ij}$. Fix an $x\in V$ and a basis $\{e_{1},\ldots,e_{m}\}$ of the tangent space $T_{x}V$. This defines a basis $\{\p_{1},\ldots, \p_{m}\}$ of $T_{\psi(x)}M_{1}$, with $\p_{i} = D\psi_{x}(e_{i})$. Viewing $\psi$ both as a map into $M_{1}$ and $\R^{d}$, we find that the metric $g$ on $M_{1}$ relative to this basis is given by:
\begin{align*}
    g_{ij}(\psi(x))
    &= g(\p_{i},\p_{j})
    = g(D\psi_{x}(e_{i}),D\psi_{x}(e_{j})  \\
    &= g_{\R^{d}}(D\psi_{x}(e_{i}),D\psi_{x}(e_{j})) \\
    &= \sum_{k=1}^{d}\frac{\p\psi_{k}}{\p x_{i}}\frac{\p\psi_{k}}{\p x_{j}}
    = [(D\psi_{x})^{t}(D\psi_{x})]_{ij};
\end{align*}
here $g_{\R^{d}}$ denotes the standard metric on $\R^{d}$. Therefore
\begin{align*}
    \psi^{*}dV_{1}
    &= \psi^{*}\paren{\det(g_{ij})^{1/2}\,dy^{1}\wedge \cdots \wedge dy^{m}} \\
    &= \det[(D\psi_{x})^{t}(D\psi_{x})]^{1/2}\,
    d(\psi^{-1}\circ \psi)^{1}\wedge\cdots \wedge d(\psi^{-1}\circ \psi)^{m}\\
    &= \det[(D\psi_{x})^{t}(D\psi_{x})]^{1/2}\,dx^1 \wedge \cdots \wedge dx^m.
\end{align*}
Thus
\begin{equation} \label{eq:lemma_change_of_var}
    \int_{U}h\,dV_{1}
    = \int_{V}(h\circ \psi)\det[(D\psi_{x})^{t}(D\psi_{x})]^{1/2}\,dx.
\end{equation}
On the other hand, since $F\circ \psi:V\to M_{2}$ is also a parametrization of $F(U)$, the same reasoning shows that
\begin{align*}
    \int_{F(U)}g\,dV_{2}
    &= \int_{V}(g\circ F\circ \psi)\det[(D(F\circ \psi)_{x})^{t}(D(F\circ \psi)_{x})]^{1/2}\,dx \\
    &= \int_{V}(g\circ F\circ \psi)\det[(D\psi_{x})^{t}(DF_{\psi(x)})^{t}(DF_{\psi(x)})(D\psi_{x})]^{1/2}\,dx.
\end{align*}
The claim then follows by taking
\begin{equation*}
    h(y)
    = g(F(y))\paren{\frac{\det[(D\psi_{\psi^{-1}(y)})^{t}(DF_{y})^{t}(DF_{y})(D\psi_{\psi^{-1}(y)})]}{\det[(D\psi_{\psi^{-1}(y)})^{t}(D\psi_{\psi^{-1}(y)})]}}^{1/2}
\end{equation*}
in formula \eqref{eq:lemma_change_of_var} above.
\end{proof}

\begin{proof} [Second Proof of Lemma \ref{lemma:delta}]
If $\eps |k|^2 + k\cdot u = 0$ and $k\neq 0$, then we can write $k = |k|\hat{k}$ to see that $\hat{k}\cdot u = -\eps |k| < 0$, implying $\hat{k}\in S_{-}^{d-1}$, and that $|k| = |\hat{k}\cdot u|/\eps$. Conversely, if $\hat{k}\in S_{-}^{d-1}$, then $\frac{|\hat{k}\cdot u|}{\eps}\hat{k}\in \wt{M}\sm \{0\}$. Therefore
\begin{equation*}
    F(\hat{k}) := \frac{|\hat{k}\cdot u|}{\eps}\hat{k} 
    = -\frac{\hat{k}\cdot u}{\eps}\hat{k}
\end{equation*}
defines a bijection from $S := \{\hat{k} \in S^{d-1} : \hat{k}\cdot u < 0\}$ to $M := \wt{M}\sm \{0\}$. Since both $F$ and its inverse (projection onto $S^{d-1}$) extend to smooth maps on open subsets of $\R^d$, we see that $F$ is a diffeomorphism.

We first focus on the case where $u = |u|e_{d}$. Applying Lemma \ref{lemma:change_of_var} to the above $F$ with $M_{1} = S$, $M_{2} = M$, $U=S$, and the parametrization 
\begin{equation*}
    \psi : \{x\in \R^{d} : |x| < 1\} \to S,
    \qquad \R^{d-1}\times \R \ni (x',x_{d}) \mapsto (x',\sqrt{|u|^{2}-|x'|^2}),
\end{equation*}
we get
\begin{equation} \label{eq:delta_proof_det}
    \frac{1}{|u|}\int_{\wt{M}}g(k)\,dV(k) \\
    = \frac{1}{|u|}\int_{S_{-}^{d-1}}g\paren{-\frac{\hat{k}\cdot u}{\eps}\hat{k}}
    \paren{\frac{\det[(D\psi_{\psi^{-1}(\hat{k})})^{t}(DF_{\hat{k}})^{t}(DF_{\hat{k}})(D\psi_{\psi^{-1}(\hat{k})})]}{\det[(D\psi_{\psi^{-1}(\hat{k})})^{t}(D\psi_{\psi^{-1}(\hat{k})})]}}^{1/2}\,d\hat{k},
\end{equation}
where we have used the fact that $M$ and $S$ differ from $\wt{M}$ and $S_{-}^{d-1}$, respectively, by sets of measure zero. Writing $\hat{k} = (k_1,\ldots,k_d)\in S_{-}^{d-1}$, we see that
\begin{equation*}
    (DF_{\hat{k}})^{t}(DF_{\hat{k}})
    = \frac{|u|^{2}}{\eps^{2}}
    \begin{pmatrix}
        k_d^{2} & 0 & \cdots & k_1 k_d \\ 
        0 & k_d^{2} & \cdots & k_2 k_d \\
        \vdots & \vdots & \ddots & \vdots \\
        k_1 k_d & k_2 k_d & \cdots & k_1^2 + \cdots + k_{d-1}^{2} + 4k_d^2
    \end{pmatrix}
    = \frac{|u|^{2}}{\eps^{2}}
    \begin{pmatrix}
        k_d^{2} & 0 & \cdots & k_1 k_d \\ 
        0 & k_d^{2} & \cdots & k_2 k_d \\
        \vdots & \vdots & \ddots & \vdots \\
        k_1 k_d & k_2 k_d & \cdots & 1 + 3k_d^2
    \end{pmatrix}
\end{equation*}
and
\begin{equation*}
    D\psi_{\psi^{-1}(\hat{k})}
    = D\psi_{(k_1,\ldots,k_{d-1})}
    = \begin{pmatrix}
        1 &  \cdots & 0 \\
        \vdots  & \ddots & \vdots \\
        0 &  \cdots & 1 \\
        -k_1/k_d & \cdots & -k_{d-1}/k_d
    \end{pmatrix},
\end{equation*}
so
\begin{align*}
    (D\psi_{\psi^{-1}(\hat{k})})^{t}(DF_{\hat{k}})^{t}(DF_{\hat{k}})(D\psi_{\psi^{-1}(\hat{k})})]
    &= \frac{|u|^2}{\eps^2}\brac{k_{d}^2 I_{d-1} + \paren{\frac{1}{k_d^2} - 2}\hat{k}'\otimes \hat{k}'} \\
    (D\psi_{\psi^{-1}(\hat{k})})^{t}(D\psi_{\psi^{-1}(\hat{k})}) &= I_{d-1} + \frac{1}{k_{d}^2}\hat{k}'\otimes \hat{k}',
\end{align*}
where $I_{d-1}$ is the $(d-1)\times (d-1)$ identity matrix and $\hat{k}' = (k_1,\ldots, k_{d-1})$. We can thus use the matrix determinant identity
\begin{equation*}
    \det (A + v\otimes w) = (1 + v\cdot A^{-1}w)\det(A)
\end{equation*}
and the fact that $\hat{k}'\cdot \hat{k}' = 1-k_d^2$ to conclude that
\begin{align*}
    \pfrac{\det[(D\psi_{\psi^{-1}(\hat{k})})^{t}(DF_{\hat{k}})^{t}(DF_{\hat{k}})(D\psi_{\psi^{-1}(\hat{k})})] }{\det[(D\psi_{\psi^{-1}(\hat{k})})^{t}(D\psi_{\psi^{-1}(\hat{k})})]}^{1/2}
    &= \pfrac{|u|}{\eps}^{d-1}\pfrac{x_d^{2(d-3)}}{1/x_d^2}^{1/2} = \pfrac{|u|}{\eps}^{d-1}|x_d|^{d-2} \\
    &= \pfrac{|u|}{\eps}^{d-1}\pfrac{|\hat{k}\cdot u|}{|u|}^{d-2}
    = \frac{|u|}{\eps^{d-1}}|\hat{k}\cdot u|^{d-2},
\end{align*}
which, when inserted into \eqref{eq:delta_proof_det}, results in the desired identity \eqref{eq:delta_id_goal}. This proves the claim for $u$ of the form $|u|e_d$, and the general case follows from the rotational symmetry of \eqref{eq:delta}.
\end{proof}

\section{Proof of Functional Estimate } \label{subsect:proof_of_functional_estimate}

In this section we provide a proof of Lemma \ref{lem:QUU_well-defined} in Subsection \ref{subsect:functional_estimate}, which is used to show the well-definedness of $Q_{qB}^{\alpha,\eps}(f)$.

\begin{proof}[Proof of Lemma \ref{lem:QUU_well-defined}]
Using Lemma \ref{lemma:delta} and the inequality $|f(v)|\le\langle v\rangle^{-N}\|f\|_{L^\infty_N}$, we can rewrite and bound $Q_{qB}^{\alpha,
\eps}(f)$ as
\begin{align*}
    \int_{\R^3}\!|Q_{qB}^{\alpha,\eps}(f)|\,dv
    &\le
    \frac{C}{\eps^4}\int_{\R^3}\!\int_{\R^3}\!\int_{S^2_-}
    |\hat k\cdot u|
    \left[
        \hat V\left(\frac{|\hat k\cdot u|}{\eps}\right)^2
        +
        \hat V\left( \frac{\sqrt{|u|^2 - |\hat k\cdot u|}}{\eps}\right)^2
    \right]
    \\
    &\qquad\qquad\qquad\cdot
    \left( 
        \langle v'\rangle^{-N}\langle{v_1'}\rangle^{-N}+\langle{v}\rangle^{-N}\langle{v_1}\rangle^{-N}
    \right)\|f\|_{L^\infty_N}^2(1+\alpha\|f\|_{L^\infty})\,d \hat k\, d v_1\, dv
    \\
    &\le
    \frac{C}{\eps^4}\|f\|_{L^\infty_N}^2(1+\alpha\|f\|_{L^\infty})
    \\
    &\qquad\qquad\cdot
    \int_{\R^3}\!\int_{\R^3}\!\int_{S^2_-}|\hat k\cdot u|\left[
        \hat V\left(\frac{|\hat k\cdot u|}{\eps}\right)^2
        +
        \hat V\left( \frac{\sqrt{|u|^2 - |\hat k\cdot u|}}{\eps}\right)^2
    \right]\langle v\rangle^{-N}\langle v_1\rangle^{-N}\,d\hat k\,d v_1\,d v,
\end{align*}
where we have used that
\begin{align*}
    \langle v\rangle^2
    \langle v_1\rangle^2
    =
    1+|v|^2+|v_1|^2+|v|^2|v_1|^2
    \le
    1+|v'|^2+|v_1'|^2+(|v'|^2+|v_1'|^2)^2
    \le 3\langle v'\rangle^2
    \langle v_1'\rangle^2
\end{align*}
to bound the weights $ \langle v'\rangle^{-N}\langle{v_1'}\rangle^{-N}\lesssim_N\langle{v}\rangle^{-N}\langle{v_1}\rangle^{-N}$. By computing the angular integral
\begin{align*}
    &\frac{1}{\eps^4}\int_{S_-^2}|\hat k\cdot u|\left[
        \hat V\left(\frac{|\hat k\cdot u|}{\eps}\right)^2
        +
        \hat V\left( \frac{\sqrt{|u|^2 - |\hat k\cdot u|}}{\eps}\right)^2
    \right]
    \,d\hat k
    \\
    &\qquad\qquad
    =\frac{2\pi}{\eps^3}\int_0^1\! \gamma\lambda\Big[\hat V(\gamma\lambda)^2 + \hat V\left(\gamma\sqrt{1-\lambda^2}\right)^2\Big]\,d\lambda
    \qquad\qquad\qquad\lambda:=|\hat k\cdot\hat u|,\quad \gamma:=|u|/\eps
    \\
    &\qquad\qquad
    =\frac{4\pi}{\eps^3}\int_0^1\! \gamma\lambda\hat V(\gamma\lambda)^2\,d\lambda
    \\
    &\qquad\qquad
    =\frac{4\pi}{\eps^3\gamma}\int_0^\gamma\!\mu\hat V(\mu)^2\,d\mu
    \le \frac{C}{\eps^2|u|}
    \int_0^\infty\!\mu\hat V(\mu)^2\,d\mu,
\end{align*}
where in the third line we have used that the change of variables $\lambda\mapsto\sqrt{1-\lambda^2}$ keeps the measure $\lambda\,d\lambda$ invariant, and combining with our previous bound, we get
\begin{align*}
\int_{\R^3}\!|Q_{qB}^{\alpha,\eps}(f)|\,dv
&\le
\frac{C}{\eps^2}\|f\|_{L^\infty_N}^2(1+\alpha\|f\|_{L^\infty})\left(
\int_0^\infty\!\mu\hat V(\mu)^2\,d\mu
\right)
\int_{\R^3}\!\int_{\R^3}\! \frac{\langle v\rangle^{-N}\langle v_1\rangle^{-N}}{|v-v_1|}d v_1\,d v
\\
&
\le\frac{C_N}{\eps^2}\|f\|_{L^\infty_N}^2(1+\alpha\|f\|_{L^\infty})
\int_0^\infty\!\mu\hat V(\mu)^2\,d\mu
,
\end{align*}
as desired.
\end{proof}

\section{Taylor Expansion} \label{subsect:taylor}

\begin{lemma} \label{lemma:taylor}
If $g\in W^{3,\infty}(\R^{3})$, then for all $v,v'\in \R^{3}$,
\begin{equation} \label{eq:taylor_3}
    g(v') - g(v)
    = (v'-v)\cdot \nabla_{v}g
    + \frac{1}{2}\Tr\paren{(v'-v)\otimes (v'-v)\,D^{2}g(v)} + |v'-v|^{3}r_{3}(v,v'),
\end{equation}
where $\|r_{3}\|_{L^\infty(\R^3\times\R^3)} \lesssim \|g\|_{W^{3,\infty}}$. If $g\in W^{2,\infty}(\R^{3})$, then for all $v,v'\in \R^{3}$,
\begin{equation} \label{eq:taylor_2}
    g(v') - g(v)
    = (v'-v)\cdot \nabla_{v}f
    + |v'-v|^{2}\rho_{2}(v,v'),
\end{equation}
where $\|\rho_2\|_{L^\infty(\R^3\times \R^3)} \lesssim \|g\|_{W^{2,\infty}}$. If $g\in W^{1,\infty}(\R^3)$, then for all $v,v'\in \R^3$, 
\begin{equation} \label{eq:taylor_1}
    g(v') - g(v) = |v'-v|\varrho_{1}(v,v'),
\end{equation}
where $\|\varrho_1\|_{L^{\infty}(\R^3\times \R^3)}\lesssim \|g\|_{W^{1,\infty}}$.
\end{lemma}
\begin{proof}
To prove \eqref{eq:taylor_3}, we first assume $g\in C^\infty$ and take a second order Taylor expansion of $g$ around $v$ to get:
\begin{align*}
    g(v')
    &= \sum_{|\alpha|\leq 2}\frac{1}{\alpha!}\p^{\alpha}g(v)(v'-v)^{\alpha} + R(v,v') \\
    &= g(v) + \nabla g(v) \cdot (v'-v)
    + \frac{1}{2}\sum_{1\leq i,j\leq 3}\frac{\p^{2}u}{\p x_{i}\p x_{j}}(v)(v'_{i}-v_{i})(v'_{j} -v_{j}) + R(v,v') \\
    &= g(v) + \nabla g(v) \cdot (v'-v) + \frac{1}{2}\operatorname{Tr}\sparen{(v' - v)\otimes (v'-v)\,D^{2}g(v)} + R(v,v'),
\end{align*}
where
\begin{align*}
    R(v,v')
    = \sum_{|\alpha|=3}\frac{3}{\alpha!}\paren{\int_{0}^{1}(1-s)^{2}\p^{\alpha}g\sparen{s(v'-v)}\,ds}(v'-v)^{\alpha}.
\end{align*}
If $v'\neq v$, we can rewrite this as 
\begin{align*}
    R(v') = |v'-v|^{3}\sum_{|\alpha|=3}\frac{3}{\alpha!}\paren{\int_{0}^{1}(1-s)^{2}\p^{\alpha}g\sparen{s(v'-v)}\,ds}\pfrac{v'-v}{|v'-v|}^{\alpha}.
\end{align*}
Thus, this Taylor expansion gives \eqref{eq:taylor_3} with
\begin{align}
    \label{eq:taylor_3_remainder}
    |r_{3}(v,v')|
    \leq \sum_{|\alpha|=3}3\|\p^{\alpha}g\|_{L^{\infty}}\paren{\int_{0}^{1}(1-s)^{2}\,ds}\abs{\frac{v'-v}{|v'-v|}}^{\alpha} 
    = \sum_{|\alpha|=3}\|\p^{\alpha}g\|_{L^{\infty}}
    = \|g\|_{W^{3,\infty}},
\end{align}
which proves the identity \eqref{eq:taylor_3} for $g\in C^\infty$. To treat the case $g\in W^{3,\infty}$, we define $g_\eps=g*\eta_\eps$ where $\eta_\eps=\eta(\cdot/\eps)\in C^\infty_c(\mathbb{R}^3)$ is an approximate identity. Since $g_\eps\in C^\infty$, the estimate \eqref{eq:taylor_3} holds for $g_\eps$, and since $g_\eps\to g$ in $C^2$ and $\|D^3g_\eps\|_{L^\infty}\le\|D^3g\|_{L^\infty}$ for all $\eps>0$, the remainder estimate \eqref{eq:taylor_3_remainder} still holds.
The proof of the other identities \eqref{eq:taylor_2} and \eqref{eq:taylor_1} follow similarly.
\end{proof}

\printbibliography

\end{document}